\documentclass{amsart}

\usepackage{graphicx}
\usepackage[font=small,labelfont=bf]{caption}
\usepackage{epstopdf}
\usepackage{enumitem}
\usepackage{amsthm}
\usepackage{verbatim}
\usepackage{tikz}
\usetikzlibrary{arrows}
\usetikzlibrary{decorations.markings}
\usepackage{yhmath}
\usepackage{relsize}
\usepackage{calc}
\usepackage{mathrsfs}

\newtheorem{theorem}{Theorem}[section]
\newtheorem{lemma}[theorem]{Lemma}
\newtheorem{corollary}[theorem]{Corollary}

\theoremstyle{definition}

\theoremstyle{remark}
\newtheorem{remark}[theorem]{Remark}

\numberwithin{equation}{section}



\begin{document}

\title{2-Systems of Arcs on Spheres with Prescribed Endpoints}

\author{Sami Douba}
\address{Department of Mathematics and Statistics, McGill University, Montreal, Quebec, Canada}
\email{sami.douba@mail.mcgill.ca}
\thanks{The author was supported in part by the grant \#346300 for IMPAN from the Simons Foundation and the matching 2015-2019 Polish MNiSW fund.}

\begin{abstract}
Let $S$ be an $n$-punctured sphere, with $n \geq 3$. We prove that $\binom{n}{3}$ is the maximum size of a family of pairwise non-homotopic  simple arcs on $S$ joining a fixed pair of distinct punctures of $S$ and pairwise intersecting at most twice. On the way, we show that a square annular diagram $A$ has a corner on each of its boundary paths if $A$ contains at least one square and the dual curves of $A$ are simple arcs joining the boundary paths of $A$ and pairwise intersecting at most once.
\end{abstract}

\maketitle

\section{Introduction}

A $k$-{\it system} $\mathcal{A}$ of arcs on a punctured surface $S$ is a collection of essential simple arcs on $S$ such that no two arcs of $\mathcal{A}$ are homotopic or intersect more than $k$ times. We begin with the following observation.

\begin{figure}
\scalebox{.45}{
\includegraphics[trim=0 520 0 0, clip]{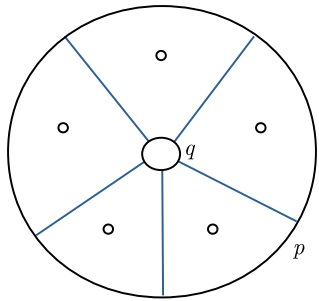}
}
\captionof{figure}{A maximum-size $0$-system joining distinct punctures $p,q$ of the $7$-punctured sphere.}
\label{fig:zerosystem}
\end{figure}

\begin{remark}\label{kequalszero}
If $k=0$, $S$ is an $n$-punctured sphere with $n \geq 3$, and the arcs of $\mathcal{A}$ all join a fixed pair of distinct punctures $p,q$ of $S$, then $|\mathcal{A}| \leq n-2$. To see this, fix a complete hyperbolic metric on $S$ of area $2\pi(n-2)$ and realize the arcs of $\mathcal{A}$ as geodesics on $S$. Cutting $S$ along $\mathcal{A}$, we obtain a collection of hyperbolic punctured strips. Since the boundary of each strip consists of two arcs of $\mathcal{A}$, and since each arc of $\mathcal{A}$ appears twice as a boundary component of some strip, we count precisely $|\mathcal{A}|$ strips. The bound on $|\mathcal{A}|$ now follows from the fact that each of these strips has area at least $2\pi$. Moreover, this bound is tight since we can easily devise a $0$-system on $S$ whose complement consists entirely of once-punctured strips (see Figure~\ref{fig:zerosystem}). An area argument also shows that the maximum size of $\mathcal{A}$ is $2n-5$ if we assume instead that $p=q$.
\end{remark}

Problems involving bounding the size of a $k$-system of arcs, of which Remark~\ref{kequalszero} serves as a trivial example, originated in similar problems for curves. Juvan, Malni\v{c}, and Mohar introduced the term ``$k$-system" and showed that the maximum size $N(k, \Sigma)$ of a $k$-system of essential simple closed curves on a fixed compact surface $\Sigma$ is finite \cite{juvan1996systems}. Independently, Farb and Leininger inquired about $N(k, g) := N(k, \Sigma)$ for $\Sigma$ closed and oriented of genus $g$ and $k=1$. In response, Malestein, Rivin, and Theran provided an upper bound exponential in $g$, and showed that $N(1,2) = 12$ \cite{malestein2014topological}. Also for $k=1$, Przytycki produced an upper bound on the order of $g^3$ \cite[Theorem~1.4]{przytycki2015arcs}; since then, tighter bounds on $N(1,g)$ have been found by Aougab, Biringer, and Gaster \cite{aougab2017packing}, and more recently, Greene \cite{greene2018curves}. Moreover, Przytycki provided an upper bound on $N(k,g)$ for arbitrary $k$ that grows like $g^{k^2+k+1}$ \cite[Cor.~1.6]{przytycki2015arcs}. This bound was subsequently improved by Greene to one that grows like $g^{k+1}\log g$ \cite{greene2018curvesii}. 

In \cite{przytycki2015arcs}, so as to prove the aforementioned results about $k$-systems of curves, Przytycki first proved stronger results about $k$-systems of arcs -- for example, that the maximum size of a $k$-system of arcs on a punctured surface $S$ of Euler characteristic $\chi < 0$ (where distinct arcs are not required to have the same endpoints) grows like $|\chi|^{k+1}$. In the same article, Przytycki proved the following:

\begin{theorem}\cite[Theorem~1.7]{przytycki2015arcs}\label{przytycki}
Let $p,q$ be punctures of an $n$-punctured sphere $S$, where $n \geq 3$. The maximum size of a $1$-system $\mathcal{A}$ of arcs on $S$ joining $p$ and $q$ is $\binom{n-1}{2}$. 
\end{theorem}

Note that $p$ and $q$ are not assumed to be distinct in the statement of Theorem~\ref{przytycki}. More recently, Bar-Natan showed that for $S,p,q$ as in Theorem~\ref{przytycki}, if $p=q$ then the maximum size of a $2$-system of arcs on $S$ joining $p$ and $q$ is $\binom{n}{3}$ \cite{bar2017arcs}. The main result of this article is that Bar-Natan's maximum holds for $p,q$ distinct:

\begin{theorem}\label{maintheorem}
Let $p,q$ be distinct punctures of an $n$-punctured sphere $S$, where $n \geq 3$. The maximum size of a $2$-system $\mathcal{A}$ of arcs on $S$ joining $p$ and $q$ is $\binom{n}{3}$. 
\end{theorem}

\subsection*{Organization} In Section \ref{example}, we provide an example of a $2$-system of size $\binom{n}{3}$ joining a fixed pair of distinct punctures of an $n$-punctured sphere for $n \geq 3$. The remaining sections are concerned with proving that $\binom{n}{3}$ is an upper bound on the size of such a $2$-system $\mathcal{A}$. This is proved by induction on $n$; we prove that the number of arcs of $\mathcal{A}$ that become homotopic after forgetting a puncture of $S$ is not too large. This, in turn, is proved by induction, via the following lemma:

\begin{lemma}\label{relation}
Let $S$ be an $n$-punctured sphere, and let $p,q,r$ be distinct punctures of $S$. Let $\mathcal{P}, \mathcal{Q}$ be $1$-systems of arcs starting at $r$ and ending at $p,q$, respectively, so that no arc of $\mathcal{P}$ intersects an arc of $\mathcal{Q}$. Suppose $\mathcal{R} \subset \mathcal{P} \times \mathcal{Q}$ satisfies the following:
\begin{enumerate}[label=(\roman*)]
\item there is at most one point of intersection between any two pairs of arcs in $\mathcal{R}$;
\item the cyclic order around $r$ of any two intersecting pairs $(\alpha, \beta), (\alpha', \beta') \in \mathcal{R}$ with $\alpha \neq \alpha', \beta \neq \beta'$ is alternating.
\end{enumerate}
Then $|\mathcal{R}| \leq \binom{n-1}{2}$. 
\end{lemma}

Section \ref{proofoflemma} is devoted to the proof of Lemma~\ref{relation}. The inductive step again involves forgetting a puncture $s$ of $S$, but this time, we choose $s$ with care in order to control the subsequent behavior of the arcs of $\mathcal{P}$ and $\mathcal{Q}$. More precisely, we require that $s$ be $p$-{\it isolated} (see Section \ref{definitions} for definitions). 

To guarantee that a puncture with this property exists, we take a detour into annular square diagrams. A $k$-{\it system annular diagram} $A$ is an annular diagram whose dual curves constitute a $k$-system of arcs joining the boundary paths of $A$. Such a diagram arises as the dual square complex to a $k$-system $\mathcal{A}$ on a punctured sphere with distinct prescribed endpoints. In Section \ref{diagrams}, we prove the following:

\begin{theorem}\label{corner}
Let $A$ be a $1$-system annular diagram. Then either $A$ is a cycle, or $A$ has a corner on each of its boundary paths.
\end{theorem}

Roughly speaking, a corner corresponds to an isolated puncture. Note that Theorem~\ref{corner} does not hold for $k=2$; see Figure~\ref{fig:diagrams} (bottom left) for a counterexample, suggested by Przytycki. 

\subsection*{Acknowledgements} I thank my supervisor Piotr Przytycki for his crucial guidance and remarkable patience throughout the course of this article's preparation. I also thank Daniel Wise for providing helpful resources.

\section{Definitions}\label{definitions}

\subsection{Arc systems} A {\it puncture} is a topological end of a space $S$ obtained from a connected, oriented, compact surface $\Sigma$ by removing finitely many points $p_1, \ldots, p_n$ from $\Sigma$. Note that the punctures of $S$ are in bijection with $p_1, \ldots, p_n$, and that we allow punctures on the boundary of $\Sigma$. If $p_1, \ldots p_n$ are taken from the interior of $\Sigma$, we call $S$ an $n$-{\it punctured} $\Sigma$. 

An {\it arc} on $S$ is a proper map $\alpha: (0,1) \rightarrow S$. A proper map induces a map between ends of topological spaces; in this sense, $\alpha$ ``maps" each endpoint of $(0,1)$ to a puncture $p$ of $S$. We call $p$ an {\it end} of $\alpha$. If $p,q$ are ends of $\alpha$, we say $\alpha$ {\it starts} at $p$ and {\it ends} at $q$, or that $\alpha$ {\it joins} $p$ and $q$. A {\it segment} of $\alpha$ is the restriction of $\alpha$ to some positive-length subinterval of $(0,1)$. 

An arc $\alpha$ is {\it simple} if it is an embedding, in which case we identify $\alpha$ and its segments with their images in $S$. If $J$ is a subinterval of $(0,1)$ with endpoints $t_1, t_2$ and $\alpha$ is a simple arc mapping $t_i$ to $x_i$ for $i=1,2$, we denote the segment $\alpha \bigr|_J$ by $(x_1x_2)_{\alpha}$. If $R \subset S$ is a subset and $p$ is an end of $\alpha$ corresponding to an endpoint $t_0=0,1$ of $(0,1)$, we say the {\it $p$-end of $\alpha$ lies in $R$} if $\alpha^{-1}(R)$ is a neighbourhood of $t_0$ in $(0,1)$.

A {\it homotopy} between arcs $\alpha_1$ and $\alpha_2$ is a proper map $(0,1) \times [0,1] \rightarrow S$ whose restrictions to $(0,1) \times \{0\}$ and $(0,1) \times \{1\}$ are $\alpha_1$ and $\alpha_2$, respectively. In particular, a homotopy preserves ends. If $r$ is a puncture of $S$, we say that two arcs on $S$ are {\it $r$-homotopic} if they are homotopic on the surface $\bar{S}$ obtained from $S$ by forgetting $r$. Two arcs are in {\it minimal position} if the number of their intersection points cannot be decreased by a homotopy. Note that if a pair of arcs have a point of intersection that is not transversal, then they are not in minimal position. An arc $\alpha$ is {\it essential} if it cannot be homotoped into a puncture, in the sense that there is no proper map $(0,1) \times [0,1) \rightarrow S$ whose restriction to $(0,1) \times \{0\}$ is $\alpha$. Unless otherwise stated, all arcs in the article are simple and essential, and all intersections between arcs are transversal. Note that an arc joining distinct punctures of a punctured surface is automatically essential.

\begin{figure}
\scalebox{.5}{
\includegraphics[trim=0 550 0 0, clip]{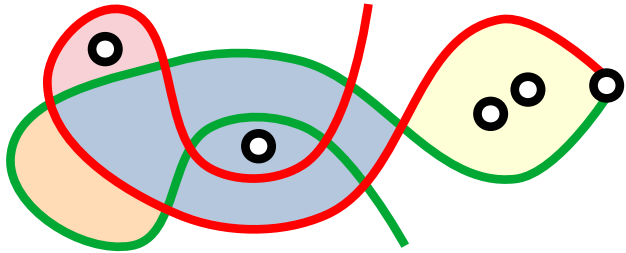}
}
\captionof{figure}{Regions formed by arcs. The yellow region is a half-bigon. The pink and orange regions are bigons. The orange bigon is empty.}
\label{fig:region}
\end{figure}

Let $R$ be a closed disc with at most 2 punctures on its boundary and possibly with punctures in its interior. A {\it region} between arcs $\alpha_1, \alpha_2$ on $S$ is a properly embedded $R \subset S$ such that $\partial R$ is a union of exactly two segments $\sigma_1, \sigma_2$, where $\sigma_i$ is a segment of $\alpha_i$ for $i=1,2$ (see Figure~\ref{fig:region}). We say that $\alpha_1, \alpha_2$ (or, more specifically, $\sigma_1, \sigma_2$) {\it form} or {\it bound} $R$. If $R$ has no punctures in its interior, we say $R$ is {\it empty}. If $R$ has exactly 0 (respectively, 1, 2) punctures on its boundary and $R \cap (\alpha_1 \cup \alpha_2) = \partial R$, we call $R$ a {\it bigon} (respectively, {\it half-bigon}, {\it strip}). We say $R$ is {\it adjacent} to a puncture $p$ of $S$ if $p$ lies on the boundary of $R$. If $p,s$ are distinct punctures of $S$ and $\mathcal{A}$ is a collection of arcs on $S$ with $s$  contained in a half-bigon or strip $H$ adjacent to $p$ formed by a pair of arcs of $\mathcal{A}$ such that $H$ is a component of $S - \bigcup \mathcal{A}$, we say that $s$ is $p$-{\it isolated by} $\mathcal{A}$. 

We will make frequent use of the following lemma.

\begin{lemma}[The bigon criterion]\label{bigoncriterion} \cite[Proposition~1.7]{farb_margalit_2012}
Two intersecting arcs on a punctured surface are in minimal position if and only if they form no empty regions.
\end{lemma}

Since an empty region bounded by intersecting arcs must contain an empty bigon or half-bigon, we immediately obtain the following corollary.

\begin{corollary}\label{corollarybigon}
Two intersecting arcs on a punctured surface are in minimal position if and only if they form no empty bigons or half-bigons.  
\end{corollary}

A {\it $k$-system} of arcs on a punctured surface $S$ is a collection $\mathcal{A}$ of essential simple arcs on $S$ such that no two arcs of $\mathcal{A}$ are homotopic or have more than $k$ points of intersection. We will mainly consider the case where $S$ is a sphere punctured at least thrice and $\mathcal{A}$ is a 2-system of arcs joining a fixed pair of distinct punctures $p,q$ of $S$. Note that for any two arcs $\alpha_1, \alpha_2$ of such a collection $\mathcal{A}$, a region bounded by $\alpha_1, \alpha_2$ that contains neither $p$ nor $q$ must be a bigon, a half-bigon, or a strip.

If a punctured surface $S$ has Euler characteristic $\chi < 0$, then $S$ admits a complete hyperbolic metric of area $2\pi |\chi|$. Under such a metric, the homotopy class of any arc contains a unique geodesic representative, and any two distinct geodesic arcs are in minimal position. Thus, for the purposes of determining the size of a $k$-system $\mathcal{A}$ of arcs on $S$, we may assume that $\mathcal{A}$ consists of geodesics.

\subsection{Combinatorial complexes} A map $X \rightarrow Y$ between CW complexes $X,Y$ is {\it combinatorial} if its restriction to each open cell of $X$ is a homeomorphism onto an open cell of $Y$. A CW complex $X$ is {\it combinatorial} if the attaching map of each  cell in $X$ is combinatorial for some subdivision of the sphere. A cell of dimension 0 is a {\it vertex} and a cell of dimension 1 is an {\it edge}. The {\it degree} of a vertex $v$ of $X$ is the number of edges in $X$ incident to $v$, with loops counted twice. 

 \subsection{Square complexes} An {$n$-cube} is a copy of $[-1,1]^n$. A {\it square complex} $X$ is a combinatorial complex whose cells are $n$-cubes with $n \leq 2$; that is, $X$ is a combinatorial 2-complex each of whose 2-cells is attached via a combinatorial map from a 4-cycle into the 1-skeleton of $X$. The cells of $X$ are called {\it cubes}, and its 2-cells are called {\it squares}. 
 
 A {\it midcube} is a subspace of a cube $[-1,1]^n$ obtained by restricting one coordinate to $0$. Let $U$ be a new square complex whose cells are midcubes of $X$ and whose attaching maps are restrictions of attaching maps in $X$ to midcubes. A {\it dual curve} $\alpha$ of a cube $c$ in $X$ is a connected component of $U$ containing a midcube of $c$. Note that if $c$ is a square, then it has exactly two dual curves. If $c$ is an edge, we say $\alpha$ is {\it dual} to $c$. We call $\alpha$ an {\it arc} if it is homeomorphic to an interval (possibly of length 0). There is a natural immersion $\alpha \rightarrow X$; if this map is an embedding, we say $\alpha$ is {\it simple}. In this case, we identify $\alpha$ with its image in $X$.

\subsection{Annular diagrams} 
\begin{figure}
\scalebox{.25}{
\includegraphics[trim=0 350 0 10, clip]{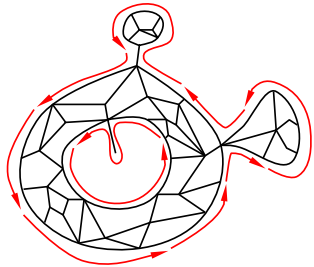}
}
\captionof{figure}{An annular diagram. The boundary paths are indicated in red.}
\label{fig:annulardiagram}
\end{figure}

An {\it annular diagram} $A$ is a finite combinatorial cell decomposition of $S^2$ minus two disjoint open 2-cells (see Figure~\ref{fig:annulardiagram}). The attaching map of each of these 2-cells is a {\it boundary path} of $A$. 

\begin{figure}
\scalebox{0.55}{
\includegraphics[trim=0 330 0 0, clip]{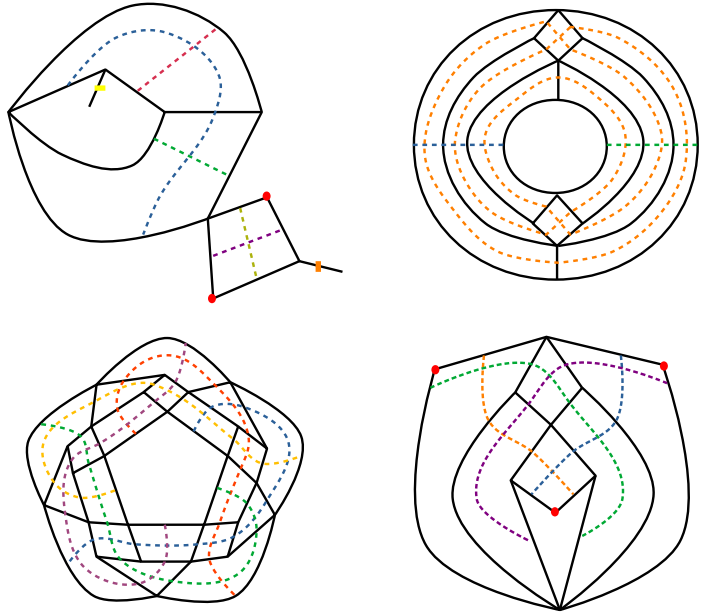}
}
\captionof{figure}{Square annular diagrams. Dual curves are dashed and colored. Corners are colored red. The bottom left (respectively, bottom right) diagram is a $2$-system (respectively, $1$-system) annular diagram. Note that the square at which the blue and orange dual curves meet in the bottom right diagram is a cornsquare with an outerpath on each boundary path, while the square at which the green and purple dual curves meet is not a cornsquare, even though the latter two dual curves are dual to consecutive edges on each boundary path.} 
\label{fig:diagrams}
\end{figure}

We call $A$ a {\it square annular diagram}, or simply a {\it diagram}, if it is also a square complex (see Figure~\ref{fig:diagrams}). A {\it corner} on a boundary path $P$ of a diagram $A$ is a vertex $v$ on $P$ of degree 2 that is contained in some square of $A$. 

Let $c$ be a square of a diagram $A$ with boundary path $P$, and let $x$ be the center of $c$. Suppose the dual curves $\alpha, \beta$ of $c$ are dual to consecutive edges $a,b$ on $P$ with shared vertex $v$. Let $\gamma$ be the loop obtained from the subarcs of $\alpha, \beta$ joining $x$ and $P$ and the half-edges of $a,b$ containing $v$. If $\gamma$ is homotopic in $A$ to a constant path, then we call $c$ a {\it cornsquare} with {\it outerpath} $ab$. 

\begin{figure}
\scalebox{.3}{
\includegraphics[trim=0 570 0 0, clip]{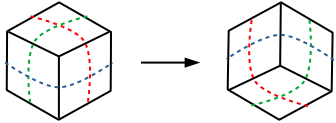}
}
\captionof{figure}{A hexagon move and its effect on dual curves.}
\label{fig:hexagon}
\end{figure}

A {\it hexagon move} on a diagram $A$ is the replacement of three squares forming a subdivided hexagon by an alternate three squares forming a subdivided hexagon (see Figure \ref{fig:hexagon}). A hexagon move can be visualized as a benign ``sliding" operation on one of the dual curves of $A$, so that if $A'$ is obtained from $A$ by a hexagon move, there is a natural correspondence between the dual curves of $A$ and those of $A'$. Note that the number of squares of $A$ is preserved under hexagon moves.

A square annular diagram $A$ is a {\it $k$-system annular diagram} if its dual curves are simple arcs joining the boundary paths of $A$ and pairwise intersecting at most $k$ times in $A$. Note that the number of intersections between any pair of dual curves of $A$ is preserved under a hexagon move. Thus, if $A'$ is obtained from a $k$-system annular diagram $A$ by a hexagon move, then $A'$ is also a $k$-system annular diagram.

\section{A 2-system of maximum size}\label{example}

We provide an example of a 2-system of arcs of size $\binom{n}{3}$ joining a fixed pair of distinct punctures of an $n$-punctured sphere $S$. This collection was independently discovered by Assaf Bar-Natan.

We think of $S$ as $\mathbb{R}^2$ punctured at $p = (-1,0)$ and at the points $r_i = (i - \frac{1}{2}, 0)$ for $i = 1, \ldots, n-2$. We construct a 2-system $\mathcal{A}$ joining $p$ and the puncture $q$ at infinity. 

Let $\alpha_{< -1}$ be the arc given by the ray $\{(x,0) : \> x < -1\}$. For $a, b, c \in \{0, 1, \ldots, n-2\}$ with $a < b <c$ or $0 < a < b = c = n-2$, let $\alpha_{abc}$ be the graph of the polynomial function $f_{abc}: (-1, \infty) \rightarrow \mathbb{R}$ given by $x \mapsto (x+1)(x-a)(x-b)(x-c)$. 

Note that for distinct triples $(a,b,c), (a', b',c')$, the difference $f_{abc}-f_{a'b'c'}$ is a cubic polynomial, one of whose roots is $-1$. Thus, the $\alpha_{abc}$ pairwise intersect at most twice. Furthermore, the $\alpha_{abc}$ are pairwise non-homotopic \cite[proof~of~Lemma~4.2]{bar2017arcs}, and $\alpha_{<-1}$ is not homotopic to any of the $\alpha_{abc}$ since the complement of $\alpha_{<-1} \cup \alpha_{abc}$ is a pair of punctured strips. 

\begin{figure}
 \begin{tikzpicture}[scale=1]

\draw [
    decoration={markings,mark=at position 1 with {\arrow[scale=2]{>}}},
    postaction={decorate},
    shorten >=0.4pt
    ]
[->] (-2,0)--(4,0); 
\draw [
    decoration={markings,mark=at position 1 with {\arrow[scale=2]{>}}},
    postaction={decorate},
    shorten >=0.4pt
    ]
    [->] (0,-3)--(0,5);
    
\begin{scope}
\clip (-2,-3) rectangle (4,5);

\draw[thick, green](0.5,0)--(1.5,0);
\draw [thick, red, domain=-1:4 , samples=500] plot (\x , {0.1*(\x+1)*\x*(\x-1)*(\x-2)});
\draw [thick, pink, domain=-1:4 , samples=500] plot (\x , {0.1*(\x+1)*\x*(\x-1)*(\x-3)});
\draw [thick, orange, domain=-1:4 , samples=500] plot (\x , {0.1*(\x+1)*\x*(\x-2)*(\x-3)});
\draw [thick, brown, domain=-1:4 , samples=500] plot (\x , {0.1*(\x+1)*(\x-1)*(\x-2)*(\x-3)});
\draw [thick, magenta, domain=-1:4 , samples=500] plot (\x , {0.1*(\x+1)*(\x-1)*(\x-3)^2});
\draw [thick, gray, domain=-1:4 , samples=500] plot (\x , {0.1*(\x+1)*(\x-2)*(\x-3)^2});
\draw[thick, blue](-1,0)--(-2,0);
\draw[thick, violet](-1,0)--(-1,2.1);
\draw[thick, violet](-1,2.1)--(1.75,2.1);
\draw[thick, violet](1.75,2.1)--(1.75,-2.2);
\draw[thick, violet](1.75,-2.2)--(-1.5,-2.2);
\draw[thick, violet](-1.5,-2.2)--(-1.5,2.6);
\draw[thick, violet](-1.5,2.6)--(2.25,2.6);
\draw[thick, violet](2.25,2.6)--(2.25,-3);

\draw [fill=white] (-1,0) circle [radius=0.05];
\draw [fill= black] (0,0) circle [radius=0.03];
\draw [fill= white] (0.5,0) circle [radius=0.05];
\draw [fill= black] (1,0) circle [radius=0.03];
\draw [fill= white] (1.5,0) circle [radius=0.05];
\draw [fill= black] (2,0) circle [radius=0.03];
\draw [fill= white] (2.5,0) circle [radius=0.05];
\draw [fill= black] (3,0) circle [radius=0.03];
\node at(-1.1,-0.6)[label=$p$]{};
\end{scope} 
\end{tikzpicture}
\captionof{figure}{The arcs $\alpha_{abc}$ on the 5-punctured sphere, together with arc $\alpha_{<-1}$, drawn in blue, arc $\alpha_{23}$, drawn in violet, and arc $\gamma_1$, drawn in green.}
\label{fig:maxfamily}
\end{figure}
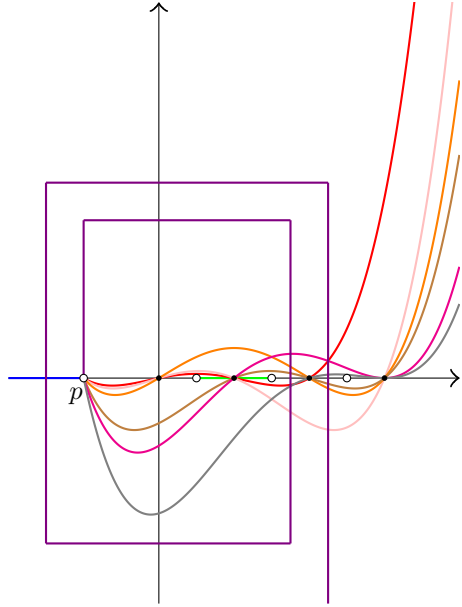

Now fix $M > 0$ such that $M > |f_{abc}(x)|$ for all $x \in (-1, n-2]$. For each $i,j \in \{1, \ldots, n-2\}$ with $i < j$, let $\alpha_{ij}$ be the union of the following horizontal and vertical segments: the segment joining $(-1,0)$ and $(-1, M)$, the segment joining $(-1,M)$ and $(i - \frac{1}{2} + \frac{1}{4}, M)$, the segment joining $(i - \frac{1}{2} + \frac{1}{4}, M)$ and $(i - \frac{1}{2} + \frac{1}{4}, -M)$, the segment joining $(i - \frac{1}{2} + \frac{1}{4}, -M)$ and $(-2, -M)$, the segment joining $(-2, -M)$ and $(-2, M+1)$, the segment joining $(-2, M+1)$ and $(j - \frac{1}{2} - \frac{1}{4}, M+1)$, and the vertical ray travelling down from $(j - \frac{1}{2} - \frac{1}{4}, M+1)$. Note that each $\alpha_{ij}$ intersects $\alpha_{<-1}$ exactly once and each $\alpha_{abc}$ exactly twice (see Figure \ref{fig:maxfamily}). Furthermore, each $\alpha_{ij}$ is in minimal position with $\alpha_{<-1}$ by Corollary~\ref{corollarybigon}; since $\alpha_{<-1}$ is disjoint from the $\alpha_{abc}$, this shows that none of the $\alpha_{ij}$ is homotopic to any of the $\alpha_{abc}$. 

We claim that the $\alpha_{ij}$ are pairwise non-homotopic. Indeed, for $k \in \{1, \ldots, n-3\}$, let $\gamma_k$ be the horizontal arc joining the punctures at $x = k- \frac{1}{2}$ and $x = k + \frac{1}{2}$, and note that $\alpha_{ij}$ and $\gamma_k$ are in minimal position by Corollary~\ref{corollarybigon}. Since no two of the $\alpha_{ij}$ share the same number of intersection points with each of the $\gamma_k$, the $\alpha_{ij}$ must be pairwise non-homotopic.

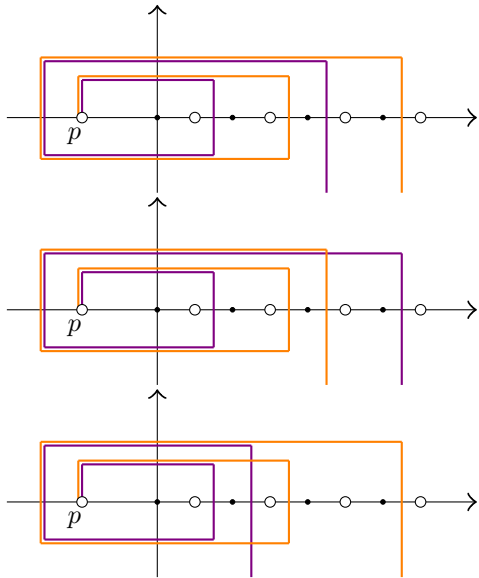
\begin{figure}
 \begin{tikzpicture}[scale=1]
\draw [
    decoration={markings,mark=at position 1 with {\arrow[scale=2]{>}}},
    postaction={decorate},
    shorten >=0.4pt
    ]
[->] (-2,0)--(4.25,0); 
\draw [
    decoration={markings,mark=at position 1 with {\arrow[scale=2]{>}}},
    postaction={decorate},
    shorten >=0.4pt
    ]
    [->] (0,-1)--(0,1.5);
    
\begin{scope}

\draw[thick, violet](-1,0)--(-1,0.5);
\draw[thick, violet](-1,0.5)--(0.75,0.5);
\draw[thick, violet](0.75,0.5)--(0.75,-0.5);
\draw[thick, violet](0.75,-0.5)--(-1.5,-0.5);
\draw[thick, violet](-1.5,-0.5)--(-1.5,0.75);
\draw[thick, violet](-1.5,0.75)--(2.25,0.75);
\draw[thick, violet](2.25,0.75)--(2.25,-1);
\draw[thick, orange](-1.05,0)--(-1.05,0.55);
\draw[thick, orange](-1.05,0.55)--(1.75,0.55);
\draw[thick, orange](1.75,0.55)--(1.75,-0.55);
\draw[thick, orange](1.75,-0.55)--(-1.55,-0.55);
\draw[thick, orange](-1.55,-0.55)--(-1.55,0.8);
\draw[thick, orange](-1.55,0.8)--(3.25,0.8);
\draw[thick, orange](3.25,0.8)--(3.25,-1);

\draw [fill=white] (-1,0) circle [radius=0.07];
\draw [fill= black] (0,0) circle [radius=0.03];
\draw [fill= white] (0.5,0) circle [radius=0.07];
\draw [fill= black] (1,0) circle [radius=0.03];
\draw [fill= white] (1.5,0) circle [radius=0.07];
\draw [fill= black] (2,0) circle [radius=0.03];
\draw [fill= white] (2.5,0) circle [radius=0.07];
\draw [fill= black] (3,0) circle [radius=0.03];
\draw [fill= white] (3.5,0) circle [radius=0.07];
\node at(-1.1,-0.6)[label=$p$]{};
\end{scope} 
\end{tikzpicture}

 \begin{tikzpicture}[scale=1]
\draw [
    decoration={markings,mark=at position 1 with {\arrow[scale=2]{>}}},
    postaction={decorate},
    shorten >=0.4pt
    ]
[->] (-2,0)--(4.25,0); 
\draw [
    decoration={markings,mark=at position 1 with {\arrow[scale=2]{>}}},
    postaction={decorate},
    shorten >=0.4pt
    ]
    [->] (0,-1)--(0,1.5);
    
\begin{scope}

\draw[thick, violet](-1,0)--(-1,0.5);
\draw[thick, violet](-1,0.5)--(0.75,0.5);
\draw[thick, violet](0.75,0.5)--(0.75,-0.5);
\draw[thick, violet](0.75,-0.5)--(-1.5,-0.5);
\draw[thick, violet](-1.5,-0.5)--(-1.5,0.75);
\draw[thick, violet](-1.5,0.75)--(3.25,0.75);
\draw[thick, violet](3.25,0.75)--(3.25,-1);
\draw[thick, orange](-1.05,0)--(-1.05,0.55);
\draw[thick, orange](-1.05,0.55)--(1.75,0.55);
\draw[thick, orange](1.75,0.55)--(1.75,-0.55);
\draw[thick, orange](1.75,-0.55)--(-1.55,-0.55);
\draw[thick, orange](-1.55,-0.55)--(-1.55,0.8);
\draw[thick, orange](-1.55,0.8)--(2.25,0.8);
\draw[thick, orange](2.25,0.8)--(2.25,-1);

\draw [fill=white] (-1,0) circle [radius=0.07];
\draw [fill= black] (0,0) circle [radius=0.03];
\draw [fill= white] (0.5,0) circle [radius=0.07];
\draw [fill= black] (1,0) circle [radius=0.03];
\draw [fill= white] (1.5,0) circle [radius=0.07];
\draw [fill= black] (2,0) circle [radius=0.03];
\draw [fill= white] (2.5,0) circle [radius=0.07];
\draw [fill= black] (3,0) circle [radius=0.03];
\draw [fill= white] (3.5,0) circle [radius=0.07];
\node at(-1.1,-0.6)[label=$p$]{};
\end{scope} 
\end{tikzpicture}

 \begin{tikzpicture}[scale=1]
\draw [
    decoration={markings,mark=at position 1 with {\arrow[scale=2]{>}}},
    postaction={decorate},
    shorten >=0.4pt
    ]
[->] (-2,0)--(4.25,0); 
\draw [
    decoration={markings,mark=at position 1 with {\arrow[scale=2]{>}}},
    postaction={decorate},
    shorten >=0.4pt
    ]
    [->] (0,-1)--(0,1.5);
    
\begin{scope}

\draw[thick, violet](-1,0)--(-1,0.5);
\draw[thick, violet](-1,0.5)--(0.75,0.5);
\draw[thick, violet](0.75,0.5)--(0.75,-0.5);
\draw[thick, violet](0.75,-0.5)--(-1.5,-0.5);
\draw[thick, violet](-1.5,-0.5)--(-1.5,0.75);
\draw[thick, violet](-1.5,0.75)--(1.25,0.75);
\draw[thick, violet](1.25,0.75)--(1.25,-1);
\draw[thick, orange](-1.05,0)--(-1.05,0.55);
\draw[thick, orange](-1.05,0.55)--(1.75,0.55);
\draw[thick, orange](1.75,0.55)--(1.75,-0.55);
\draw[thick, orange](1.75,-0.55)--(-1.55,-0.55);
\draw[thick, orange](-1.55,-0.55)--(-1.55,0.8);
\draw[thick, orange](-1.55,0.8)--(3.25,0.8);
\draw[thick, orange](3.25,0.8)--(3.25,-1);

\draw [fill=white] (-1,0) circle [radius=0.07];
\draw [fill= black] (0,0) circle [radius=0.03];
\draw [fill= white] (0.5,0) circle [radius=0.07];
\draw [fill= black] (1,0) circle [radius=0.03];
\draw [fill= white] (1.5,0) circle [radius=0.07];
\draw [fill= black] (2,0) circle [radius=0.03];
\draw [fill= white] (2.5,0) circle [radius=0.07];
\draw [fill= black] (3,0) circle [radius=0.03];
\draw [fill= white] (3.5,0) circle [radius=0.07];
\node at(-1.1,-0.6)[label=$p$]{};
\end{scope} 
\end{tikzpicture}

\captionof{figure}{The $\alpha_{ij}$ pairwise intersect at most twice.}
\label{fig:twistingcurves}
\end{figure}

We claim further that the $\alpha_{ij}$ pairwise intersect at most twice. Indeed, if $i,j,i',j' \in \{1, \ldots, n-2\}$ with $i \leq i'$, then the number of intersection points between $\alpha_{ij}, \alpha_{i'j'}$ is determined by the order of $j, i', j'$. If $i' < j \leq j'$, then $\alpha_{ij}$ and $\alpha_{i'j'}$ are disjoint (see Figure \ref{fig:twistingcurves}, top). If $i' < j' < j$, then $\alpha_{ij}$ and $\alpha_{i'j'}$ intersect once (see Figure \ref{fig:twistingcurves}, middle). Otherwise, $j \leq i'$, and there are two points of intersection between $\alpha_{ij}$ and $\alpha_{i',j'}$ (see Figure \ref{fig:twistingcurves}, bottom). Thus, the family $\mathcal{A}$ consisting of $\alpha_{<-1}$, the $\alpha_{abc}$, and the $\alpha_{ij}$ is a 2-system of size $1 + (n-3) + \binom{n-1}{3} + \binom{n-2}{2} = \binom{n}{3}$. 

\section{Properties of $r$-homotopic arcs intersecting at most twice}

Let $p,q,r$ be distinct punctures of a punctured sphere $S$, and let $\mathcal{A}$ be a $2$-system of arcs on $S$ joining $p$ and $q$. Let $\bar{S}$ be the surface obtained from $S$ by forgetting $r$, and for each arc $\alpha \in \mathcal{A}$, let $\bar{\alpha}$ be the homotopy class of the corresponding arc on $\bar{S}$. In order to bound the size of $\mathcal{A}$ from above, we will need to examine to what extent the map $\alpha \mapsto \bar{\alpha}$ is injective. In this section, we collect some facts about the fibers of this map. Together, the results of this section show that we can extend $\mathcal{A}$ so that the size of each fiber is $1$ larger than the number of pairs of disjoint arcs in that fiber. 

The main results of this section are Lemmas~\ref{rhomotopic}, \ref{rhomotopicbigon}, and \ref{extend}. The proofs are rather technical and may be skipped on an initial reading.

\begin{figure}
\scalebox{.4}{
\includegraphics[trim=0 250 0 0, clip]{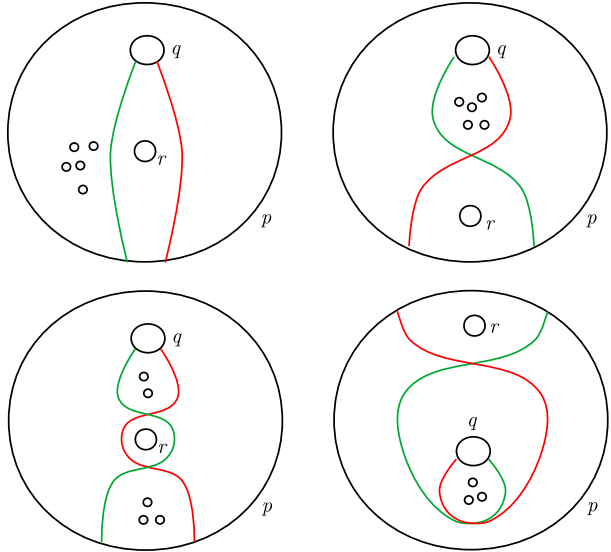}
}
\captionof{figure}{The possible configurations of a pair of $r$-homotopic arcs in minimal position and intersecting at most twice, up to relabeling $p$ and $q$.}
\label{fig:rhomotopicpairs}
\end{figure}

\begin{lemma}\label{configurations}
Let $p,q,r$ be distinct punctures of a punctured sphere $S$, and let $\alpha_1, \alpha_2$ be a pair of $r$-homotopic arcs joining $p$ and $q$ and intersecting at most twice. If the $\alpha_i$ are in minimal position, then they are in one of the configurations shown in Figure \ref{fig:rhomotopicpairs}, up to relabeling $p$ and $q$. 
\end{lemma}

\begin{figure}
\scalebox{.45}{
\includegraphics[trim=0 520 0 0, clip]{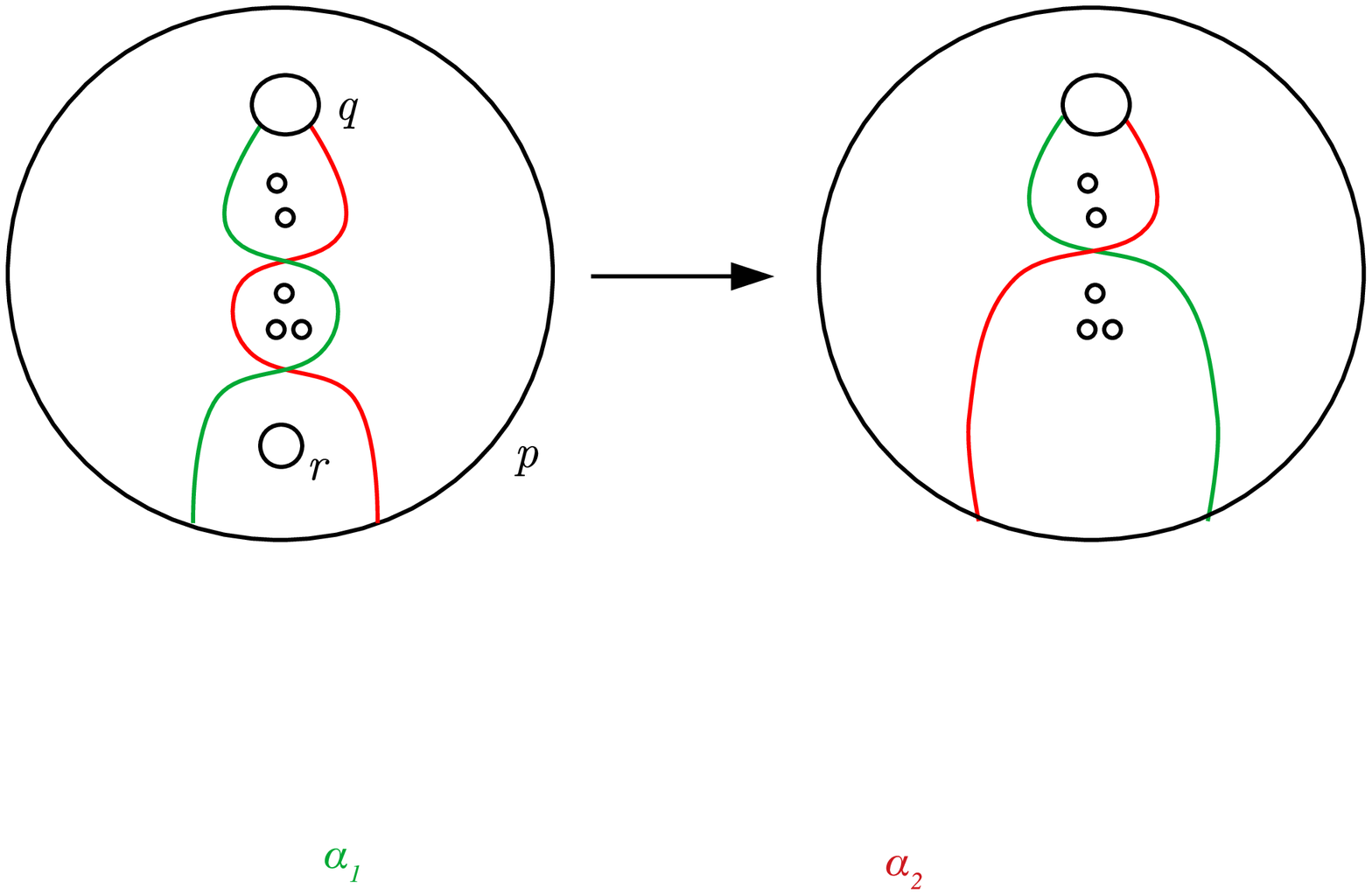}
}
\captionof{figure}{If $\alpha_1, \alpha_2$ are in minimal position, intersect exactly twice, and form a bigon that does not contain $r$, then they cannot be $r$-homotopic. The fact that the bigon and half-bigons bounded by the $\alpha_i$ prior to forgetting $r$ are punctured, and the fact that the arcs on the right are in minimal position, are consequences of Corollary~\ref{corollarybigon}.}
\label{fig:nobigonallowed}
\end{figure}

\begin{proof}
If $\alpha_1, \alpha_2$ are disjoint, then they bound a strip whose only puncture is $r$ (see Figure \ref{fig:rhomotopicpairs}, top left). Otherwise, by Corollary~\ref{corollarybigon}, $\alpha_1, \alpha_2$ bound a half-bigon or bigon $R$ whose only puncture is $r$. If $\alpha_1, \alpha_2$ intersect exactly once, then $R$ is a half-bigon and, since the $\alpha_i$ are $r$-homotopic, all punctures of $S$ distinct from $p,q,r$ lie in the other half-bigon formed by $\alpha_1, \alpha_2$ (see Figure \ref{fig:rhomotopicpairs}, top right). 

If the $\alpha_i$ intersect twice and $R$ is a half-bigon, then the $\alpha_i$ do not form a bigon, since otherwise they would not be $r$-homotopic (see Figure~\ref{fig:nobigonallowed}). Thus, in this case, the $\alpha_i$ must be as in the bottom right diagram of Figure \ref{fig:rhomotopicpairs}, and since the $\alpha_i$ are $r$-homotopic, all punctures of $S$ distinct from $p,q,r$ must lie in the other half-bigon formed by the $\alpha_i$.

Otherwise, $R$ is a bigon, and the $\alpha_i$ also bound a pair of punctured half-bigons. These half-bigons must contain all the remaining punctures of $S$ since the $\alpha_i$ are $r$-homotopic (see Figure \ref{fig:rhomotopicpairs}, bottom left).
\end{proof}

The corollary of the following lemma will be useful in the proofs of Lemmas~\ref{rhomotopic} and \ref{rhomotopicbigon}. The former tells us that, in a particular context, if we have a portion of an arc then we can trace out the remainder of that arc.

\begin{figure}
\scalebox{.4}{
\includegraphics[trim=0 520 0 0, clip]{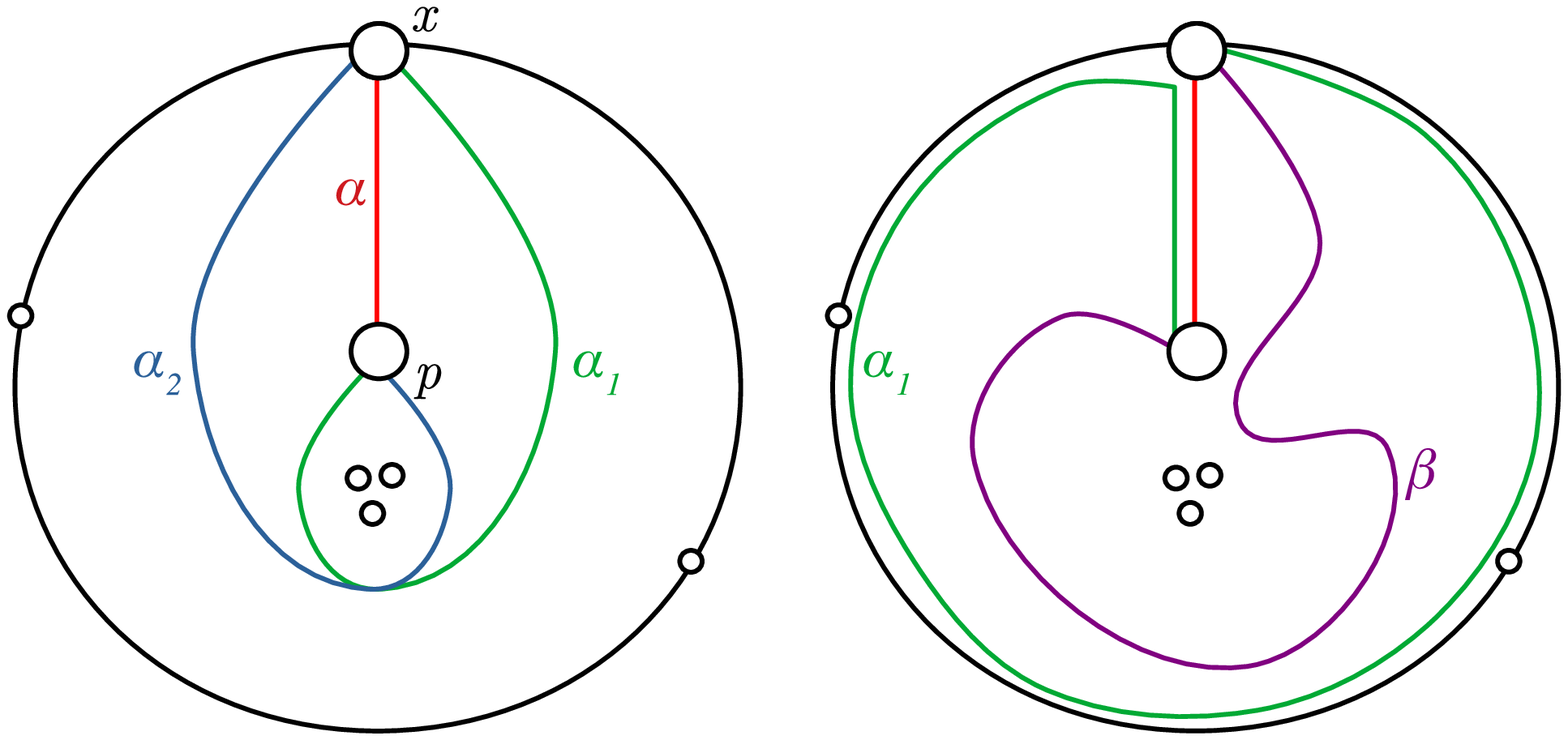}
}
\captionof{figure}{}
\label{fig:punctureddisc}
\end{figure}

\begin{lemma}\label{punctureddisc}
Let $D$ be a disc with at least $2$ punctures in its interior and at least $1$ puncture on its boundary, and let $\alpha$ be an arc joining an interior puncture $p$ of $D$ to a puncture $x$ on $\partial D$. If $\beta$ is another arc joining $p$ and $x$ such that $\alpha$ and $\beta$ bound a strip containing all interior punctures of $D$ distinct from $p$, then $\beta$ is homotopic to exactly one of the arcs $\alpha_1, \alpha_2$ shown in Figure \ref{fig:punctureddisc} (left).
\end{lemma}

\begin{proof}
Suppose that the $x$-end of $\beta$ lies to the right of $\alpha$. Then we may homotope $\alpha_1$ so that it bounds an empty strip with $\beta$, as in Figure \ref{fig:punctureddisc} (right). Thus, in this case, $\beta$ is homotopic to $\alpha_1$. Similarly, if the $x$-end of $\beta$ lies to the left of $\alpha$, then $\beta$ is homotopic to $\alpha_2$.
\end{proof}

\begin{corollary}\label{uniquecontinuation}
Let $p,q,r$ be distinct punctures of an $n$-punctured sphere $S$, with $n \geq 4$, and let $\alpha, \beta$ be $r$-homotopic arcs in minimal position joining $p,q$ and intersecting once or twice. Let $x_1, \ldots, x_m$ be the points of intersection of $\alpha, \beta$ in the order that $\beta$ traverses them as $\beta$ travels from $p$ to $q$, and set $x_0=p$, $x_{m+1}=q$. For $i=0, \ldots, m$, let $\beta_i$ be the segment of $\beta$ joining $x_i$ and $x_{i+1}$. If $m=2$, then the homotopy types of $\beta_0$ and $\beta_1$ determine that of $\beta$. If $\alpha$ and $\beta$ do not bound a bigon, then the homotopy type of $\beta_0$ determines that of $\beta$ for $m=1,2$.
\end{corollary}

\begin{proof} For $i = 0, \ldots, m$, let $\alpha_i$ be the segment of $\alpha$ joining $x_i$ and $x_{i+1}$. We puncture $S$ at $x_1, \ldots, x_m$. \\

\paragraph{\textbf{Case 1: $\alpha$ and $\beta$ intersect exactly once.}} Cutting $S$ along $\alpha_0, \beta_0$ yields two punctured strips. Let $D$ be the strip containing $q$. Note that $x_1$ is now a puncture on $\partial D$, and that $\alpha_1$ and $\beta_1$ are arcs joining $q$ and $x_1$ and bounding a strip containing all the interior punctures of $D$ distinct from $q$. Thus, by Lemma~\ref{punctureddisc}, the homotopy type of $\beta_1$ is uniquely determined, since only one of the arcs described in Lemma~\ref{punctureddisc} produces a $\beta$ that intersects $\alpha$ transversally at $x_1$ (in fact, the only other candidate homotopy class of $\beta_1$ produces a $\beta$ that is homotopic to $\alpha$). \\

\begin{figure}
\scalebox{.4}{
\includegraphics[trim=0 520 0 0, clip]{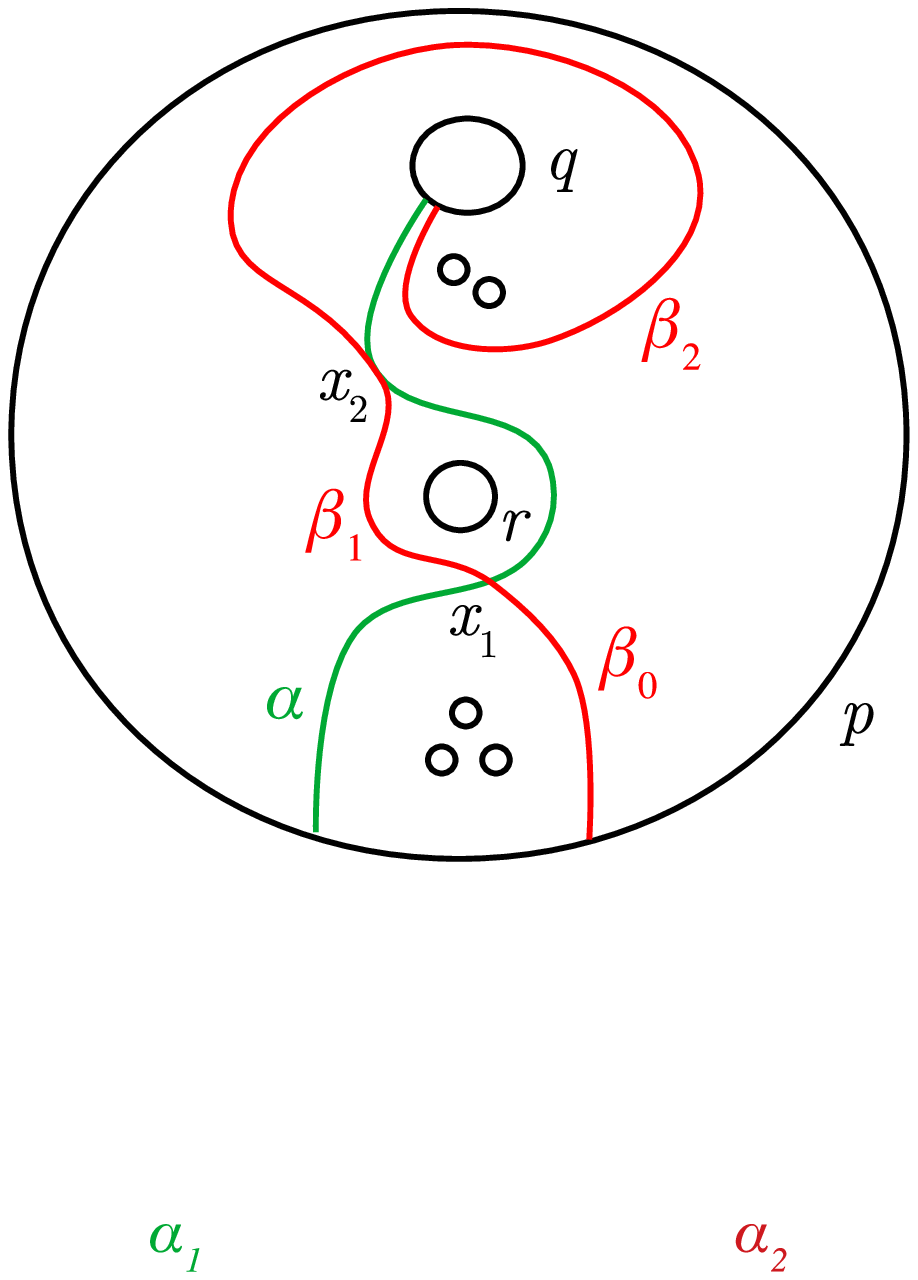}
}
\captionof{figure}{By Lemma~\ref{punctureddisc}, given segments $\beta_0$ and $\beta_1$ of $\beta$, there are at most $2$ homotopy classes of arcs joining $q$ and $x_2$ to which $\beta_2$ can belong. One such homotopy class produces a $\beta$ that intersects $\alpha$ non-transversally at $x_2$, as shown above.}
\label{fig:bigoncontinuation}
\end{figure}

\paragraph{\textbf{Case 2: $\alpha$ and $\beta$ form a bigon.}} Let $D$ be the square containing the puncture $q$ obtained by cutting $S$ along $\alpha_0, \beta_0, \alpha_1, \beta_1$. Now $x_2$ is a puncture on $\partial D$, and $\alpha_2, \beta_2$ are arcs joining $q$ and $x_2$ and bounding a strip containing all the interior punctures of $D$ distinct from $q$, so we may apply Lemma~\ref{punctureddisc} as in Case 1. Again, only one of the two homotopy classes to which $\beta_2$ must belong by Lemma~\ref{punctureddisc} produces a $\beta$ that intersects $\alpha$ transversally at $x_2$ (see Figure~\ref{fig:bigoncontinuation}).  \\

\paragraph{\textbf{Case 3: $\alpha$ and $\beta$ intersect exactly twice but do not form a bigon.}} Let $D$ be the strip containing $q$ obtained by cutting $S$ along $\alpha_0, \beta_0$. Since $x_1$ is a puncture on $\partial D$, and $\alpha_1, \beta_1$ are arcs joining $x_1, x_2$ and bounding a strip containing all the interior punctures of $D$ distinct from $x_2$, the homotopy type of $\beta_1$ is uniquely determined by Lemma~\ref{punctureddisc} as in the previous cases. Now let $D'$ be the strip containing $q$ obtained by cutting $D$ along $\alpha_1, \beta_1$. Since $x_2$ is a puncture on $\partial D'$, and $\alpha_2, \beta_2$ are arcs joining $q, x_2$ and bounding a strip containing all interior punctures of $D'$ distinct from $q$, the homotopy type of $\beta_2$ is uniquely determined by Lemma~\ref{punctureddisc}. 
\end{proof}

\begin{figure}
\scalebox{.5}{
\includegraphics[trim=0 520 0 0, clip]{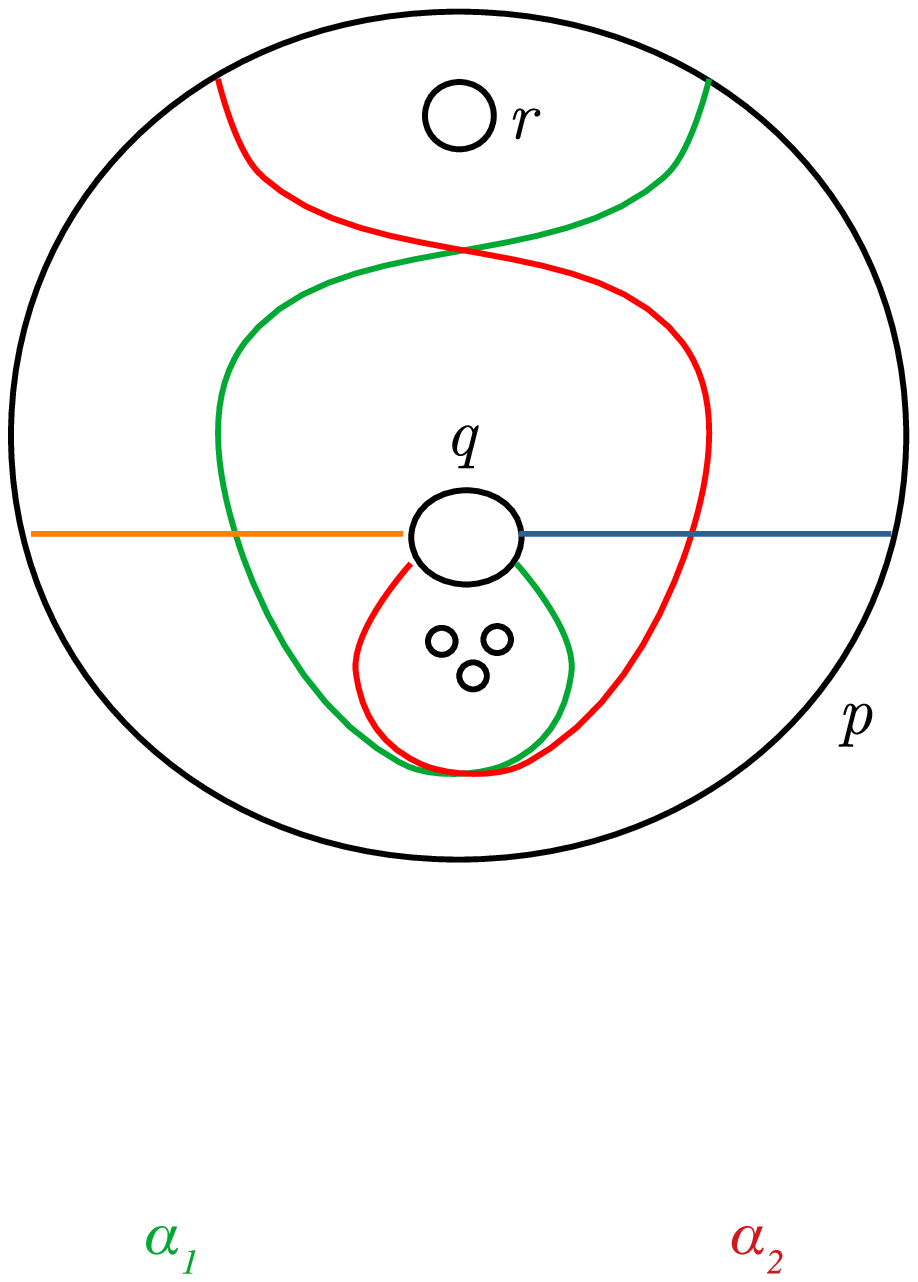}
}
\captionof{figure}{}
\label{fig:halfbigonmaximal}
\end{figure}

\begin{lemma}\label{rhomotopic}
Let $p,q,r$ be distinct punctures of an $n$-punctured sphere $S$, with $n \geq 4$. Let $\mathcal{A}_r$ be a maximal $2$-system of $r$-homotopic arcs on $S$ joining $p$ and $q$. If $\mathcal{A}_r$ contains intersecting arcs $\alpha_1, \alpha_2 \in \mathcal{A}_r$ in minimal position that do not form a bigon, then $\mathcal{A}_r$ is as in Figure \ref{fig:halfbigonmaximal}, up to homotopy and relabeling $p$ and $q$.
\end{lemma}

\begin{figure}
\scalebox{.4}{
\includegraphics[trim=0 520 0 0, clip]{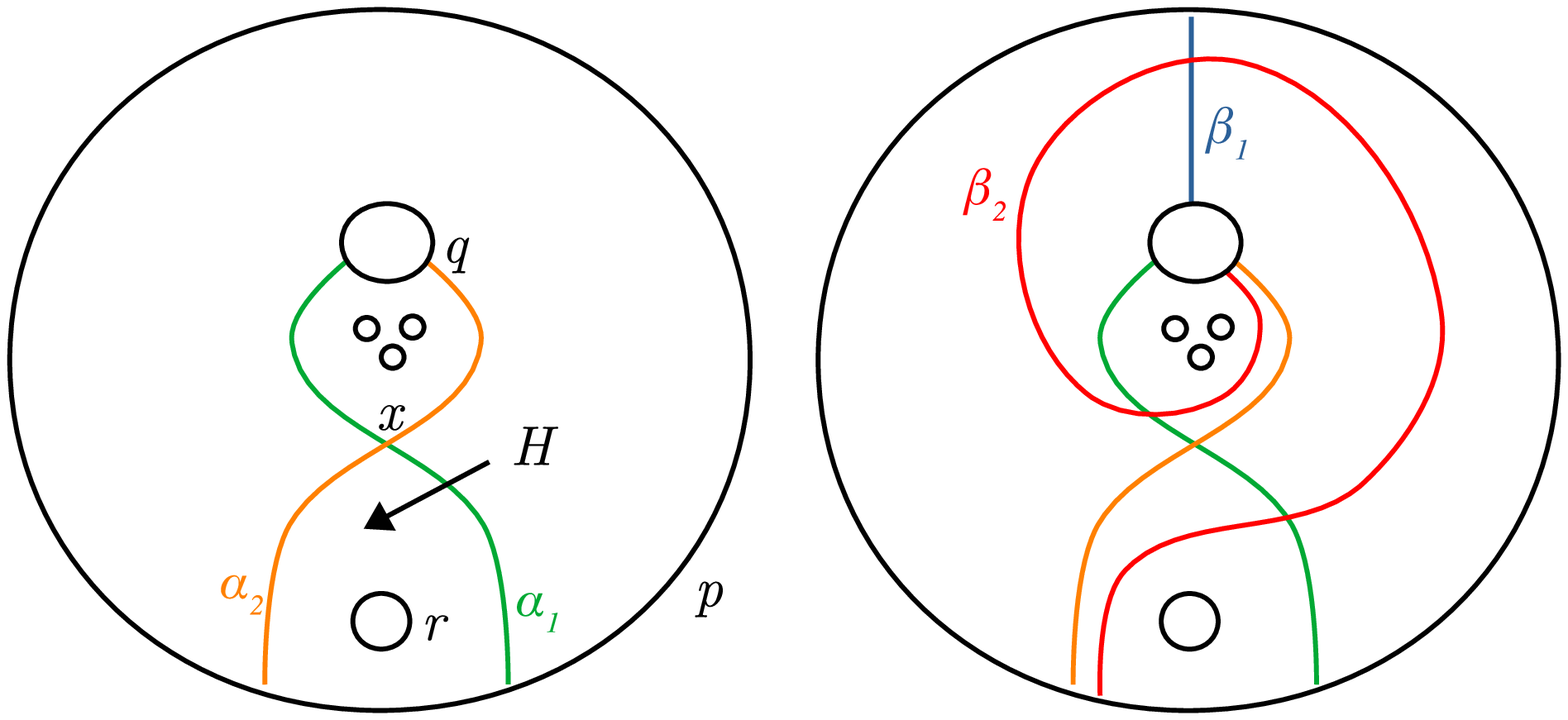}
}
\captionof{figure}{}
\label{fig:case1nobigon}
\end{figure}

\begin{figure}
\scalebox{.4}{
\includegraphics[trim=0 520 0 0, clip]{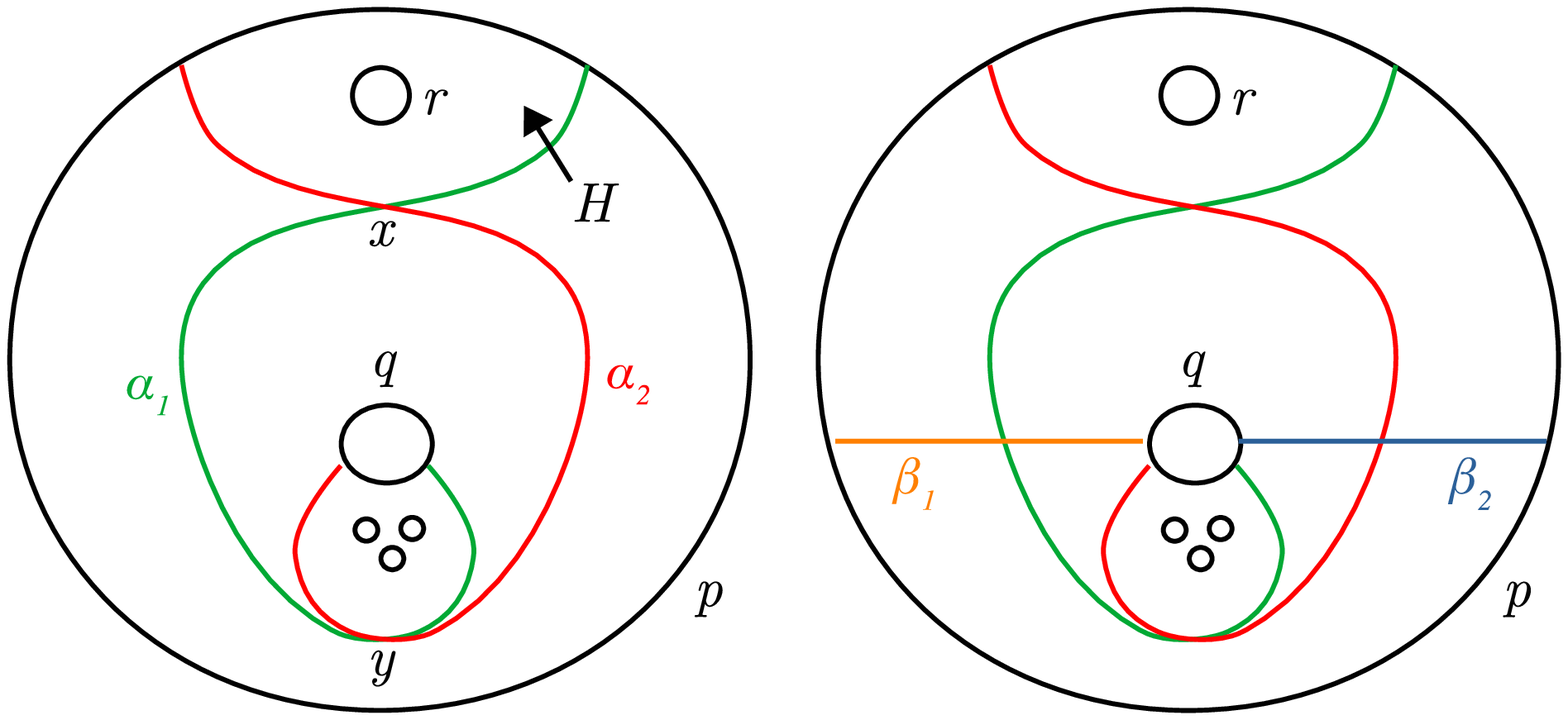}
}
\captionof{figure}{}
\label{fig:case2nobigon}
\end{figure}

\begin{proof}
Let $H$ be the half-bigon formed by the $\alpha_i$ containing $r$, and assume that $H$ is adjacent to $p$ (see Figures \ref{fig:case1nobigon}, \ref{fig:case2nobigon}, left). Let $\beta \in \mathcal{A}_r$. \\

\item \paragraph{\textbf{Case 1: The $\alpha_i$ intersect exactly once.}} Let $x$ be their unique point of intersection. If $\beta$ is disjoint from the $\alpha_i$, then $\beta$ is homotopic to arc $\beta_1$ in Figure \ref{fig:case1nobigon} (right). Now suppose $\beta$ is not disjoint from the $\alpha_i$, and let $z$ be the first point of intersection of $\beta$ and the $\alpha_i$ as $\beta$ travels from $p$ to $q$. Suppose that $z$ lies on $\alpha_1$, and that $\alpha_1, \beta$ are in minimal position. If the $p$-end of $\beta$ lies in $H$, then $\beta$ forms a half-bigon with $\alpha_1$ whose only puncture is $r$, since otherwise $\alpha_1$ and $\beta$ would form an empty half-bigon, contradicting our assumption that $\alpha_1, \beta$ are in minimal position (Corollary~\ref{corollarybigon}). Thus, by Lemma~\ref{configurations} and Corollary~\ref{uniquecontinuation}, $\beta$ is either homotopic to $\alpha_2$ or to arc $\beta_2$ in Figure \ref{fig:case1nobigon} (right). If the $p$-end of $\beta$ lies outside $H$, then $z$ cannot lie on the segment $(px)_{\alpha_1}$, since otherwise $\beta$ and $\alpha_1$ would form an empty half-bigon. Thus, $z$ lies on $(xq)_{\alpha_1}$, and so $\beta$ again forms a half-bigon with $\alpha_1$ whose only puncture is $r$. Thus, by Corollary~\ref{uniquecontinuation}, $\beta$ is again homotopic to one of $\alpha_2, \beta_2$. Note that, by the above, $\mathcal{A}_r$ cannot contain an additional arc $\beta'$ in minimal position with $\alpha_2$ and intersecting $\alpha_2$ first as $\beta'$ travels from $p$ to $q$. This is because the reflection of $\beta_2$ across the vertical diameter in Figure \ref{fig:case1nobigon} (right) intersects $\beta_2$ thrice. 

\item \paragraph{\textbf{Case 2: The $\alpha_i$ intersect exactly twice.}} Let $x,y$ be the points of intersection of the $\alpha_i$ in the order that $\alpha_1$ traverses them as it travels from $p$ to $q$. In this case, $\beta$ intersects at least one of the $\alpha_i$ since $p,q$ are in distinct components of the complement of $\alpha_1 \cup \alpha_2$. Let $z$ be the first point of intersection of $\beta$ and the $\alpha_i$ as $\beta$ travels from $p$ to $q$. We suppose again that $z$ lies on $\alpha_1$, and that $\alpha_1, \beta$ are in minimal position. If the $p$-end of $\beta$ lies in $H$, then, as in Case 1, $\beta$ forms a half-bigon with $\alpha_1$ whose only puncture is $r$. Thus, by Corollary \ref{uniquecontinuation}, $\beta$ is either homotopic to $\alpha_2$ or to arc $\beta_1$ in Figure \ref{fig:case2nobigon} (right). If the $p$-end of $\beta$ lies outside $H$, then, as in Case 1, $z$ cannot lie on the segment $(px)_{\alpha_1}$. Thus, $z$ lies on $(xy)_{\alpha_1}$. But then $\beta$ again forms a half-bigon with $\alpha_1$ whose only puncture is $r$, and so $\beta$ is either homotopic to $\alpha_2$ or to $\beta_1$ as before. Similarly, if $\beta$ is in minimal position with $\alpha_2$ and intersects $\alpha_2$ first as it travels from $p$ to $q$, then $\beta$ is either homotopic to $\alpha_1$ or to arc $\beta_2$ in Figure \ref{fig:case2nobigon} (right). 
\end{proof}

\begin{figure}
\scalebox{.5}{
\includegraphics[trim=0 520 0 0, clip]{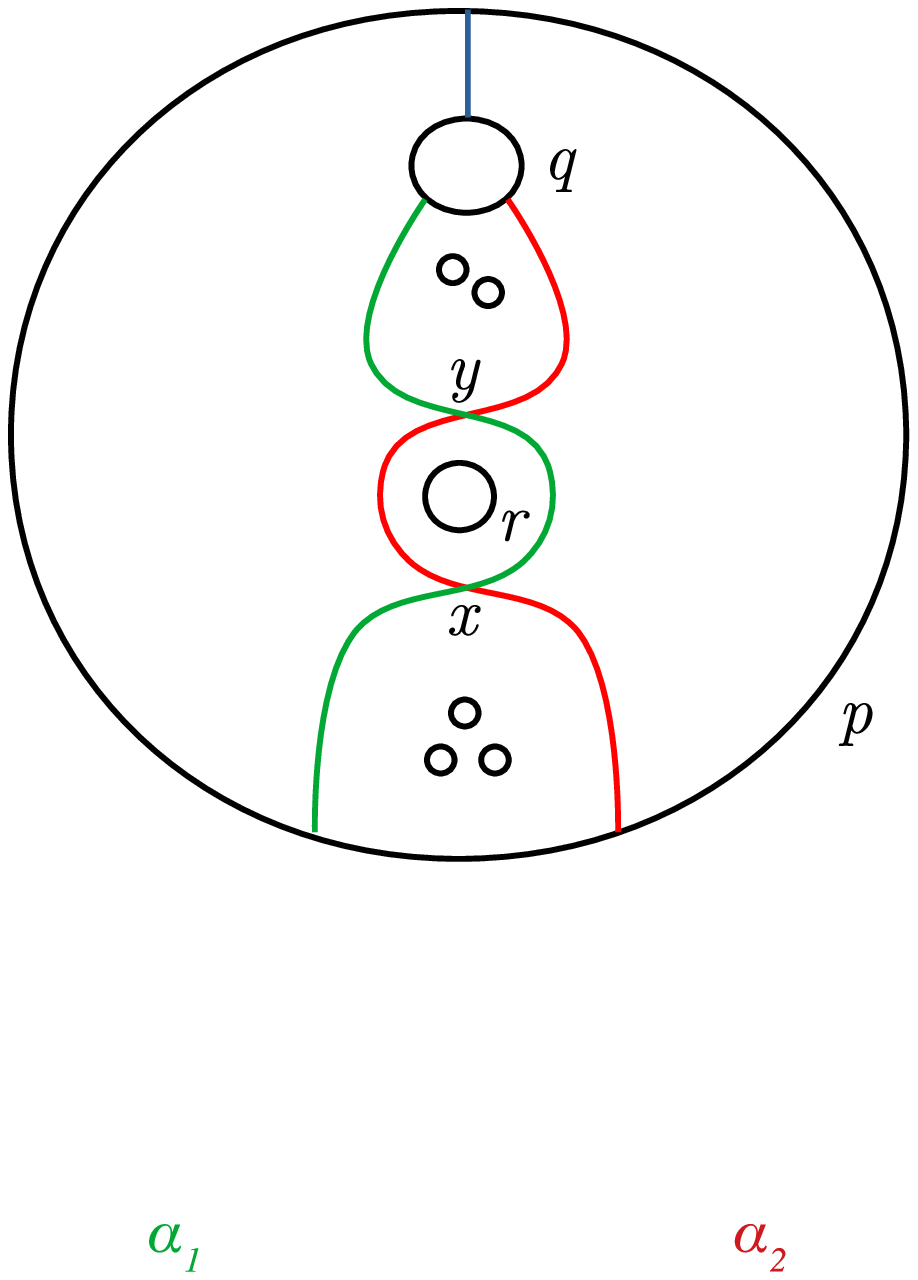}
}
\captionof{figure}{}
\label{fig:bigon}
\end{figure}

\begin{lemma}\label{rhomotopicbigon}
Let $p,q,r$ be distinct punctures of an $n$-punctured sphere $S$, with $n \geq 4$. Let $\mathcal{A}_r$ be a maximal $2$-system of $r$-homotopic arcs on $S$ joining $p$ and $q$. If $\mathcal{A}_r$ contains arcs $\alpha_1, \alpha_2 \in \mathcal{A}_r$ in minimal position that form a bigon, then $\mathcal{A}_r$ is as in Figure \ref{fig:bigon}, up to homotopy and relabeling $p$ and $q$.
\end{lemma}

\begin{figure}
\scalebox{.4}{
\includegraphics[trim=0 520 0 0, clip]{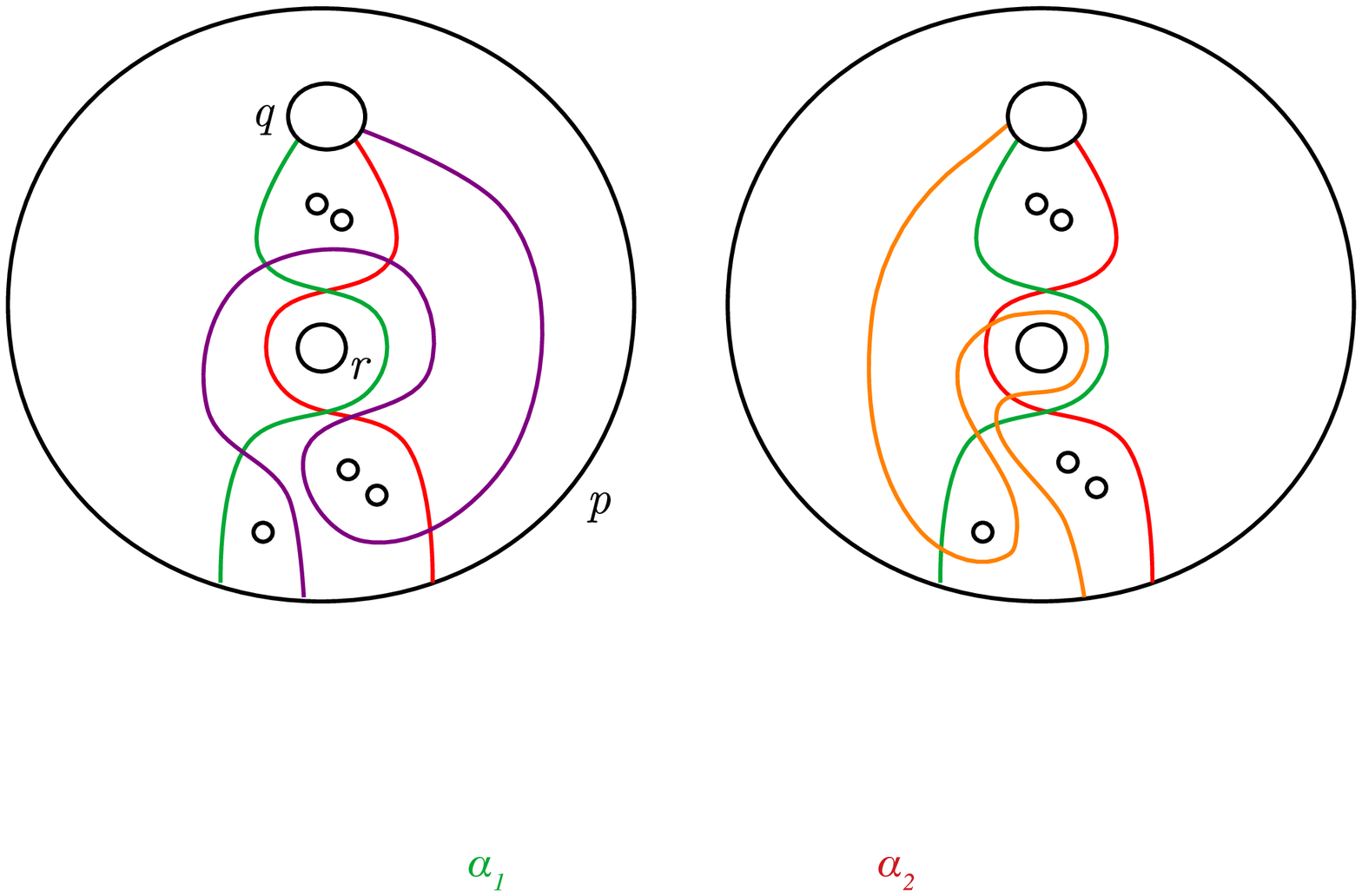}
}
\captionof{figure}{}
\label{fig:bigoncases}
\end{figure}

\begin{proof}
Let $H$ be the half-bigon adjacent to $p$ formed by the $\alpha_i$. Let $x,y$ be the points of intersection of the $\alpha_i$ in the order that $\alpha_1$ traverses them as it travels from $p$ to $q$. Let $\beta \in \mathcal{A}_r$, and assume $\alpha_1, \alpha_2, \beta$ are in minimal position. 

If $\beta$ is disjoint from the $\alpha_i$, then $\beta$ is homotopic to the blue arc in Figure \ref{fig:bigon}. Now suppose $\beta$ is not disjoint from the $\alpha_i$, and let $z_1, z_2, \ldots$ be the points of intersection of $\beta$ and the $\alpha_i$ in the order that $\beta$ traverses them as it travels from $p$ to $q$. We assume that $z_1$ lies on $\alpha_1$. Note that, by Lemma \ref{rhomotopic}, $\beta$ forms a bigon (containing only the puncture $r$) with each of the $\alpha_i$ that it intersects. \\

\item \paragraph{\textbf{Case 1: $z_1$ lies on the segment $(yq)_{\alpha_1}$.}} In this case, $\alpha_1$ and $\beta$ form a half-bigon whose only puncture is $r$. As remarked above, this is impossible.

\item \paragraph{\textbf{Case 2: $z_1$ lies on the segment $(xy)_{\alpha_1}$.}} In this case, $z_2$ does not lie on $(xy)_{\alpha_2}$. Otherwise, since $\alpha_2$ and $\beta$ are in minimal position, they would form a half-bigon whose only puncture is $r$ (as in Case 1 of Lemma~\ref{rhomotopic}), but this is impossible. Thus, $z_2$ lies on $(xy)_{\alpha_1}$, and so $\beta$ is homotopic to $\alpha_2$ by Corollary~\ref{uniquecontinuation}.

\item \paragraph{\textbf{Case 3: $z_1$ lies on the segment $(px)_{\alpha_1}$.}} In this case, since $\alpha_1, \beta$ are in minimal position, $\beta$ forms a half-bigon $H'$ with $\alpha_1$ adjacent to $p$ and containing at least one of the punctures of $H$. 

Observe that $z_2$ cannot lie on the segment $(yq)_{\alpha_2}$, since otherwise $\alpha_2$ and $\beta$ would form a half-bigon containing $r$. We also have that $z_2$ cannot lie on $(xy)_{\alpha_1}$, since otherwise $\alpha_1$ and $\beta$ would not form a bigon. Furthermore, if $z_2$ lies on $(px)_{\alpha_2}$, then so must $z_3$, since otherwise $\alpha_1$ and $\beta$ would not form a bigon. But if $z_2$ and $z_3$ both lie on $(px)_{\alpha_2}$, then $\alpha_2$ and $\beta$ form a bigon that does not contain $r$, which is impossible. 

Now, if $H'$ contains all the punctures of $H$, then $z_2$ cannot lie on $(xy)_{\alpha_2}$ since $\alpha_2, \beta$ are in minimal position, and if $z_2$ lies on $(yq)_{\alpha_1}$ then $\beta$ is homotopic to $\alpha_2$ by Corollary~\ref{uniquecontinuation}. Thus, we may assume that $H'$ contains some but not all of the punctures of $H$. 

Under this assumption, $z_2$ cannot lie on $(yq)_{\alpha_1}$, since otherwise $\beta$ would be homotopic to the purple arc Figure~\ref{fig:bigoncases} (left) by Corollary~\ref{uniquecontinuation}, and so $\beta$ would intersect $\alpha_2$ thrice. For the same reason, $z_3$ cannot lie on $(xy)_{\alpha_1}$ if $z_2$ lies on $(xy)_{\alpha_2}$. The only case left to consider is that $z_2$ and $z_3$ both lie on $(xy)_{\alpha_2}$. But then $\beta$ is homotopic to the orange arc in Figure~\ref{fig:bigoncases} (right) by Corollary~\ref{uniquecontinuation}, and so $\beta$ intersects $\alpha_1$ thrice. 
\end{proof}

\begin{figure}
\scalebox{.5}{
\includegraphics[trim=0 600 0 0, clip]{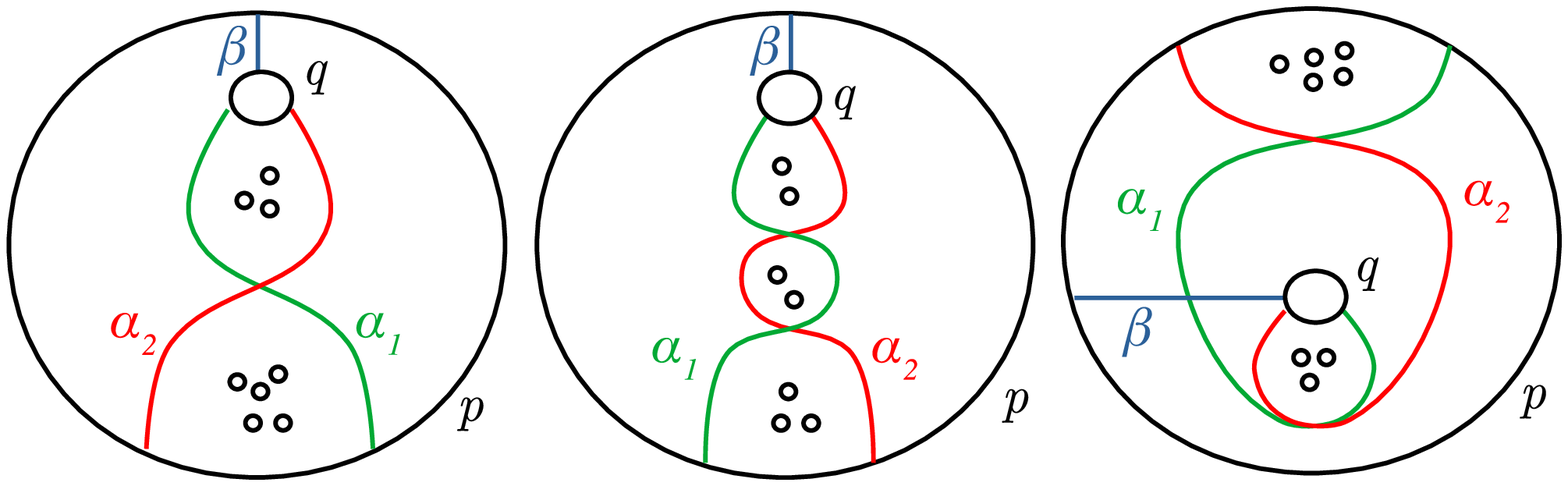}
}
\captionof{figure}{}
\label{fig:extend}
\end{figure}

\begin{lemma}\label{extend}
Let $p,q$ be distinct punctures of an $n$-punctured sphere $S$, with $n \geq 4$. Let $\alpha_1, \alpha_2, \beta$ be arcs on $S$ joining $p,q$ in one of the configurations shown in Figure~\ref{fig:extend}. Then an arc $\gamma$ joining $p,q$ that is in minimal position with $\beta$ and intersects $\beta$ at least thrice must intersect $\alpha_1$ or $\alpha_2$ at least thrice.
\end{lemma}

\begin{proof}
Set $x_0= p$, and let $x_1, x_2, x_3$ be the first 3 points of intersection of $\beta$ and $\gamma$ in the order that $\gamma$ traverses them as $\gamma$ travels from $p$ to $q$. For $i=0,1,2$, let $\beta_i = (x_ix_{i+1})_\beta, \gamma_i = (x_ix_{i+1})_\gamma$, and let $R_i$ be the region not containing $p$ bounded by $\beta_i, \gamma_i$. If $\gamma_0$ does not intersect the $\alpha_i$, then $R_0$ contains no punctures, and so $\beta, \gamma$ are not in minimal position by Lemma~\ref{bigoncriterion}, contradicting our assumption. Thus, $\gamma_0$ has at least one point of intersection with the $\alpha_i$. Similarly, $R_1, R_2$ must each contain at least one puncture of $S$, and so each of $\gamma_1, \gamma_2$ has at least 2 points of intersection with the $\alpha_i$. Thus, $\gamma$ has at least 3 points of intersection with $\alpha_1$ or $\alpha_2$. 
\end{proof}

\section{1-System annular diagrams}\label{diagrams}

In this section, we prove Theorem~\ref{corner}, which will be useful in the inductive step of the proof of Lemma~\ref{relation}.

\begin{lemma}
Let $S$ be a twice-punctured sphere and $p,q$ its punctures. Let $\mathcal{A}$ be a finite collection of simple arcs joining $p$ and $q$ and pairwise intersecting at most once. If there is a pair of intersecting arcs of $\mathcal{A}$, then there is a pair of arcs of $\mathcal{A}$ forming a half-bigon $H$ adjacent to $p$ such that no other arc of $\mathcal{A}$ has its $p$-end in $H$. 
\end{lemma}

\begin{proof}
Pick $\alpha, \beta \in \mathcal{A}$ such that $\alpha$ and $\beta$ intersect, and let $H$ be the half-bigon adjacent to $p$ formed by $\alpha, \beta$. If $H$ is as in the statement of the lemma, then we are done. Otherwise, there is an arc $\beta' \in \mathcal{A}$ whose $p$-end lies in $H$. Since the $q$-end of $\beta'$ is outside $H$, the arc $\beta'$ must intersect one of $\alpha, \beta$, say $\alpha$. We now repeat the above steps with arcs $\alpha, \beta'$. Since there are finitely many arcs in $\mathcal{A}$, this process must terminate. 
\end{proof}

The following corollary follows immediately.

\begin{corollary}\label{cornsquare}
Let $A$ be a $1$-system annular diagram, and $P$ a boundary path of $A$. If $A$ has at least one square, then $A$ has a cornsquare with outerpath on $P$.
\end{corollary}

We now proceed to the proof of Theorem~\ref{corner}. 

\begin{proof}[Proof of Theorem~\ref{corner}]
We proceed by induction on the number of squares of $A$. If $A$ has no squares, then $A$ is a cycle. Now suppose $A$ has at least one square, that there is a boundary path $P$ of $A$ without a corner, and that the theorem holds for any annular square complex with fewer squares than $A$. Since $A$ contains a square, $A$ contains a cornsquare with outerpath on $P$ by Corollary~\ref{cornsquare}. Thus, we may produce a corner on $P$ via a series of hexagon moves \cite[Figure~3.17]{wise2012riches}. Note that a single hexagon move cannot produce two corners on $P$; otherwise there would be a dual curve beginning and terminating at $P$ (see Figure \ref{fig:singlecorner}). 

\begin{figure}
\scalebox{.3}{
\includegraphics[trim=0 580 0 0, clip]{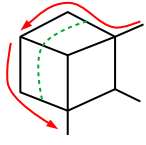}
}
\captionof{figure}{If a single hexagon move produces two corners on $P$, then there is a dual curve beginning and ending at $P$.}
\label{fig:singlecorner}
\end{figure}

We perform hexagon moves until the first corner $v$ on $P$ is produced. Note that each neighbor of $v$ has degree at least 4. Indeed, since $P$ had no corners, we had to have performed at least one hexagon move to obtain $v$, but a neighbor of $v$ of degree 3 would correspond to a corner prior to performing that move, contradicting our assumption that $v$ is the first corner produced on $P$ (see Figure~\ref{fig:vertexdegree}). Thus, by deleting $v$ as well as the two edges and the square incident to $v$, we obtain a 1-system annular diagram with one fewer square than $A$ and without any corners on one of its boundary paths, contradicting the induction hypothesis. 
\end{proof}

\begin{figure}
\scalebox{.3}{
\includegraphics[trim=0 580 0 0, clip]{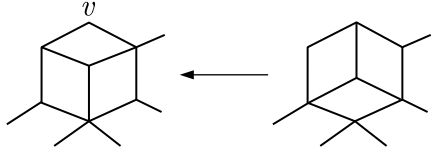}
}
\captionof{figure}{If a corner $v$ produced by a hexagon move has a neighbor of degree 3, then that neighbor had to have been a corner prior to performing that move.}
\label{fig:vertexdegree}
\end{figure}

\begin{corollary}\label{isolatedpuncture}
Let $p,q$ be punctures of a punctured sphere $S$, and let $\mathcal{A}$ be a 1-system of arcs joining $p$ and $q$ such that $|\mathcal{A}| \geq 2$ and the arcs of $\mathcal{A}$ are pairwise in minimal position. There is a puncture $s$ of $S$ distinct from $p,q$ that is $p$-isolated by $\mathcal{A}$. If $\mathcal{A}$ contains a pair of intersecting arcs, then $s$ can be chosen so that the component $H$ of $S - \bigcup \mathcal{A}$ containing $s$ is a half-bigon.
\end{corollary}

\begin{proof}
The dual square complex $A$ to $\mathcal{A}$ is an annular square complex as in Theorem~\ref{corner}. Let $P$ be the boundary path corresponding to $p$. Note that $A$ has at least 2 vertices since $|\mathcal{A}| \geq 2$, and that $A$ has at least one square if and only if $\mathcal{A}$ contains at least one pair of intersecting arcs. If $A$ is a cycle, then we may take $H$ to be the strip corresponding to any vertex of $A$. Otherwise, $A$ has a corner $v$ on $P$, and we may take $H$ to be the half-bigon corresponding to $v$. In either case, $H$ is punctured since the arcs of $\mathcal{A}$ are pairwise non-homotopic and in minimal position.
\end{proof}

\section{Proof of Lemma~\ref{relation}}\label{proofoflemma}

In this section, we prove Lemma~\ref{relation}, which essentially constitutes the inductive step in the proof of Theorem~\ref{maintheorem}. We will need the following:

\begin{lemma}\label{erdos}\cite{erdos1946sets}
A set of pairwise intersecting straight line segments between $\ell$ points on a circle in $\mathbb{R}^2$ has size at most $\ell$. 
\end{lemma}

\begin{proof}[Proof of Lemma~\ref{relation}] We fix a complete hyperbolic metric on $S$ of area $2\pi(n-2)$. We may assume that $\mathcal{P}, \mathcal{Q}$ are nonempty, and that the arcs of $\mathcal{P} \cup \mathcal{Q}$ are pairwise in minimal position. We divide the proof into steps:\\

\item \paragraph{\textbf{Step 0. The arcs of $\mathcal{P}$ (and hence the arcs of $\mathcal{Q}$) are consecutive at $r$.}}
Indeed, suppose $\alpha, \alpha' \in \mathcal{P}$ are distinct, and suppose there is an arc $\beta \in \mathcal{Q}$ whose $r$-end lies in the strip or half-bigon $H$ bounded by $\alpha, \alpha'$ and adjacent to $r$. Since $\beta$ does not intersect $\alpha, \alpha'$, the puncture $q$ must lie in $H$. Since no arc of $\mathcal{Q}$ intersects $\alpha, \alpha'$, it follows that the $r$-end of every arc of $\mathcal{Q}$ must also lie in $H$.

\item \paragraph{We fix an orientation on $S$.} This induces a cyclic order $C$ of the arcs of $\mathcal{P} \cup \mathcal{Q}$ around $r$. By Step 0, this order in turn induces a linear order $<$ on $\mathcal{P}$, where the minimum and maximum arcs of $\mathcal{P}$ are those with a successor or predecessor in $\mathcal{Q}$ under $C$. 

We proceed by induction. If $n=3$, then, up to homotopy, there is a unique arc joining $r$ to each of $p,q$, and the statement of the lemma holds. Now let $n \geq 4$, and assume the lemma holds if $S$ has fewer punctures. If $\mathcal{P}$ consists of a single arc then the lemma is trivially satisfied since $|\mathcal{Q}| \leq \binom{n-1}{2}$ by Theorem~\ref{przytycki}.  Thus, we may assume that $|\mathcal{P}| \geq 2$.

\item \paragraph{\textbf{Step 1. There is a puncture $s$ of $S$ distinct from $p,q,r$ that is $p$-isolated by $\mathcal{P} \cup \mathcal{Q}$.}} Indeed, if the arcs of $\mathcal{P}$ are pairwise disjoint, then since $|\mathcal{P}| \geq 2$, we have that $S - \mathcal{P}$ consists of at least two punctured strips adjacent to $p$ and $r$, and so we may take $s$ to be a puncture of any such strip that does not contain $q$ (see Figure~\ref{fig:step1}, left). Otherwise, by Corollary~\ref{isolatedpuncture}, there is a puncture $s$ distinct from $p,r$ that is $p$-isolated by $\mathcal{P}$ such that the component of $S - \bigcup\mathcal{P}$ containing $s$ is a half-bigon (see Figure~\ref{fig:step1}, right). In this case, $s$ is necessarily distinct from $q$ since we are assuming $\mathcal{Q}$ to be nonempty, and so there is at least one arc disjoint from the arcs of $\mathcal{P}$ joining $q$ and $r$.

\begin{figure}
\scalebox{.4}{
\includegraphics[trim=0 520 0 0, clip]{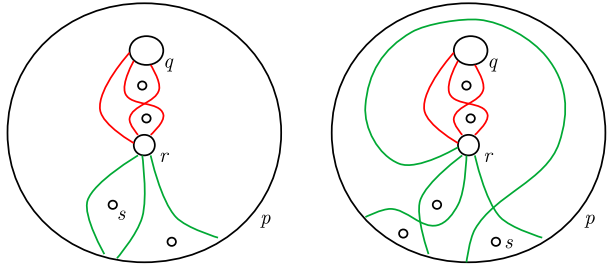}
}
\captionof{figure}{An illustration of Step 1.}
\label{fig:step1}
\end{figure}

\item \paragraph{Let $\bar{S}$ be the surface obtained from $S$ by forgetting the puncture $s$ (endowed with a complete, finite-area hyperbolic metric), and for each arc $\alpha \in \mathcal{P} \cup \mathcal{Q}$, let $\bar{\alpha}$ be the corresponding arc on $\bar{S}$.} Let $\bar{\mathcal{P}}, \bar{\mathcal{Q}}$ be the collection of all $\bar{\alpha}$ for $\alpha \in \mathcal{P}, \mathcal{Q}$, respectively. We tighten the arcs of $\bar{\mathcal{P}}\cup \bar{\mathcal{Q}}$ to geodesics, thereby identifying arcs that correspond to $s$-homotopic arcs on $S$.

The orientation on $S$ induces an orientation on $\bar{S}$. As above, this gives us a linear order $\prec$ on $\bar{\mathcal{P}}$.

\item \paragraph{\textbf{Step 2. Two distinct arcs in $\mathcal{Q}$ cannot be $s$-homotopic.}} Otherwise, they would form a strip adjacent to $q,r$ or a half-bigon adjacent to one of $q,r$ whose only puncture is $s$, which cannot happen since the component of $S - \bigcup(\mathcal{P}\cup\mathcal{Q})$ containing $s$ is adjacent to $p$. Thus, we may identify $\bar{\mathcal{Q}}$ with $\mathcal{Q}$. 

\item \paragraph{\textbf{Step 3. Arcs in $\mathcal{P}$ that are $s$-homotopic must be consecutive at $r$.}} Indeed, suppose that $\alpha, \alpha' \in \mathcal{P}$ are $s$-homotopic. If $\alpha, \alpha'$ bound a strip whose only puncture is $s$, then the $r$-end of any arc $\beta \in \mathcal{P}$ distinct from $\alpha, \alpha'$ cannot lie inside this strip, since otherwise $\beta$ and one of $\alpha, \alpha'$ would necessarily form a half-bigon adjacent to $r$ whose only puncture is $s$, contradicting the fact that $s$ is $p$-isolated. 

Otherwise, $\alpha$ and $\alpha'$ form a half-bigon adjacent to $p$ whose only puncture is $s$, and a half-bigon adjacent to $r$ containing all punctures of $S$ except $p,r,s$. Thus, any arc $\alpha''$ in $\mathcal{P}$ distinct from $\alpha, \alpha'$ with $r$-end outside the latter half-bigon must be disjoint from $\alpha,\alpha'$, in which case $\alpha, \alpha', \alpha''$ are in fact $s$-homotopic. 

\item \paragraph{\textbf{Step 4. For any $\alpha, \alpha' \in \mathcal{P}$, if $\alpha < \alpha'$ then $\bar{\alpha} \preceq \bar{\alpha'}$.}} Indeed, if $\bar{\alpha} \succ \bar{\alpha'}$, then $\alpha, \alpha'$ form a half-bigon adjacent to $r$ whose only puncture is $s$. This cannot happen since $s$ is $p$-isolated.

\item \paragraph{We now extend $\mathcal{P}$ and $\mathcal{R}$ as follows.} Suppose $\mathcal{P}$ contains intersecting $s$-homotopic arcs $\alpha, \alpha'$. Then, as discussed in Step 3, an arc $\alpha''$ disjoint from $\alpha, \alpha'$ joining $p$ and $r$ is $s$-homotopic to $\alpha, \alpha'$ and is between $\alpha, \alpha'$ with respect to the linear order on $\mathcal{P}$. If $\alpha''$ is not homotopic to an arc in $\mathcal{P}$, then we add $\alpha''$ to $\mathcal{P}$ (otherwise, name the latter arc $\alpha''$). If there is an arc $\beta \in \mathcal{Q}$ such that $(\alpha, \beta), (\alpha', \beta) \in \mathcal{R}$ but $(\alpha'', \beta) \notin \mathcal{R}$, then we add the pair $(\alpha'', \beta)$ to $\mathcal{R}$. Since $\alpha''$ is disjoint from all arcs in $\mathcal{P}\cup \mathcal{Q}$, we have not violated any of the conditions of the lemma. We do this for each pair of intersecting $s$-homotopic arcs $\alpha, \alpha'$.

Let $\bar{\mathcal{R}}$ be the image of $\mathcal{R}$ under the map $\mathcal{P} \times \mathcal{Q} \rightarrow \bar{\mathcal{P}} \times \mathcal{Q}, (\alpha, \beta) \mapsto (\bar{\alpha}, \beta)$.
It is clear that $\bar{\mathcal{R}}$ satisfies condition $(i)$, and it follows from Step 4 that $\bar{\mathcal{R}}$ satisfies $(ii)$ as well, so that $|\bar{\mathcal{R}}| \leq \binom{n-2}{2}$. It remains to show that $|\mathcal{R}| - |\bar{\mathcal{R}}| \leq n-2$. 

Let $\mathcal{I}$ be the subset of $\bar{\mathcal{R}}$ consisting of elements with more than one pre-image under the map $\mathcal{R} \rightarrow \bar{\mathcal{R}}$. 

\item \paragraph{\textbf{Step 5. If $(\bar{\alpha}, \beta), (\bar{\alpha'}, \beta') \in \mathcal{I}$ with $\beta \neq \beta'$, then $\beta, \beta'$ are disjoint.}} Indeed, let $(\alpha_1, \beta), (\alpha_2, \beta)$ be pre-images of $(\bar{\alpha}, \beta)$, and $(\alpha_1', \beta'), (\alpha_2', \beta')$ pre-images of $(\bar{\alpha'}, \beta')$ with $\alpha_1 \neq \alpha_2, \alpha_1' \neq \alpha_2'$. Suppose that $\beta$ and $\beta'$ intersect. Note that we cannot have $\{\alpha_1, \alpha_2\} = \{\alpha_1', \alpha_2'\}$, since otherwise either $(\alpha_1, \beta), (\alpha_2, \beta')$ or $(\alpha_1, \beta'), (\alpha_2, \beta)$ would be two intersecting pairs of arcs in $\mathcal{R}$ whose cyclic order around $r$ is not alternating, contradicting assumption $(ii)$. 

Now suppose $\{\alpha_1, \alpha_2\} \cap \{\alpha_1', \alpha_2'\} \neq \emptyset$. Then there are 3 distinct arcs among $\alpha_1, \alpha_2, \alpha_1', \alpha_2'$ that are $s$-homotopic. Assume without loss of generality that these arcs are $\alpha_1, \alpha_2, \alpha_1'$. Observe that $\alpha_1'$ cannot intersect $\alpha_i$ for $i=1,2$, since otherwise $(\alpha_i, \beta), (\alpha_1', \beta')$ would intersect more than once, contradicting assumption $(i)$. Thus, $\alpha_1$ intersects $\alpha_2$, but then $\alpha_1'$ lies between $\alpha_1$ and $\alpha_2$ at $r$, and so there is an $i\in \{1,2\}$ such that $(\alpha_i, \beta), (\alpha_1', \beta')$ contradict assumption $(ii)$.

We conclude that $\alpha_1, \alpha_2, \alpha_1', \alpha_2'$ are distinct; in particular, since at most 3 distinct arcs in $\mathcal{P}$ can be $s$-homotopic, we must have $\bar{\alpha} \neq \bar{\alpha'}$. Thus, by Step 3, the order of $\alpha_1, \alpha_2, \alpha_1', \alpha_2'$ at $r$ is neither alternating nor nested. By the latter, we mean that $\alpha_1, \alpha_2$ do not both lie between $\alpha_1', \alpha_2'$ in the linear order on $\mathcal{P}$, and vice versa. Since each of the pairs of arcs $\alpha_1, \alpha_2$ and $\alpha_1', \alpha_2'$ bound a strip or half-bigon containing $s$, we must have that $\alpha_i, \alpha_j'$ intersect for some $i,j \in \{1,2\}$, but then $(\alpha_i, \beta), (\alpha_j', \beta')$ intersect more than once, contradicting assumption $(i)$.

\item \paragraph{\textbf{Step 6. If $(\bar{\alpha}, \beta), (\bar{\alpha'}, \beta') \in \mathcal{I}$ with $\bar{\alpha} \neq \bar{\alpha'}, \beta \neq \beta'$, then the cyclic order of $(\bar{\alpha}, \beta), (\bar{\alpha'}, \beta')$ at $r$ is alternating.}} Indeed, suppose otherwise, and let $(\alpha_1, \beta), (\alpha_2, \beta)$ be pre-images of $(\bar{\alpha}, \beta)$, and $(\alpha_1', \beta'), (\alpha_2', \beta')$ pre-images of $(\bar{\alpha'}, \beta')$ with $\alpha_1 \neq \alpha_2, \alpha_1' \neq \alpha_2'$. Then, by Step 3, the order of $\alpha_1, \alpha_2, \alpha_1', \alpha_2'$ at $r$ is neither alternating nor nested. Thus, as in Step 5, $\alpha_i, \alpha_j'$ must intersect for some $i,j \in \{1,2\}$, but then $(\alpha_i, \beta), (\alpha_j', \beta')$ are two intersecting pairs of arcs in $\mathcal{R}$ whose cyclic order around $r$ is not alternating, contradicting assumption $(ii)$ (see Figure~\ref{fig:nestedalternating}).

\begin{figure}
\scalebox{.5}{
\includegraphics[trim=0 520 0 0, clip]{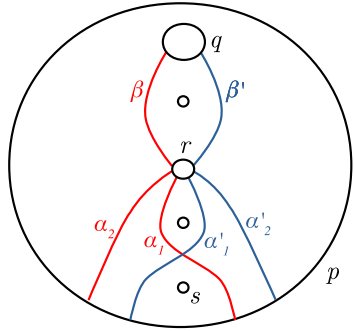}
}
\captionof{figure}{An illustration of Step 6. Here, the pairs $(\alpha_1, \beta)$ and $(\alpha_1', \beta')$ contradict assumption $(ii)$.}
\label{fig:nestedalternating}
\end{figure}

\item \paragraph{Let $\mathcal{H}_q$ be the image of $\mathcal{I}$ under the projection map $\bar{\mathcal{P}} \times \mathcal{Q} \rightarrow \mathcal{Q}$.} Let $\mathcal{H}_s$ be the collection of all geodesic arcs $a$ joining $r$ and $s$ such that $a$ is contained in a strip bounded by a pair of distinct, disjoint $s$-homotopic arcs in $\mathcal{P}$. 

 Let $\mathcal{H} = \mathcal{H}_s \cup \mathcal{H}_q$, and let $\mathcal{I}' \subset \mathcal{H}_s \times \mathcal{H}_q$ be the set of all $(a, \beta) \in \mathcal{H}_s \times \mathcal{H}_q$ such that $(\alpha, \beta) \in \mathcal{I}$ for an arc $\alpha$ bounding a strip corresponding to $a$. We extended $\mathcal{R}$ so that the map $\mathcal{R} \rightarrow \bar{\mathcal{R}}$ is injective outside a set of cardinality $|\mathcal{I}'|$. Thus, to complete the proof, it suffices to show that $|\mathcal{I}'| \leq n-2$. 

\item \paragraph{\textbf{Step 7. The complement of $\mathcal{H}$ consists of punctured strips and a single square, possibly with no punctures.}} Indeed, the arcs of $\mathcal{H}_s$ are disjoint by construction \cite[proof~of~Theorem~1.7]{przytycki2015arcs}, and the arcs of $\mathcal{H}_q$ are disjoint by Step 5, so that the complement of each of $\mathcal{H}_s, \mathcal{H}_q$ consists of punctured strips. By Step 0, the arcs in $\mathcal{H}_q$ (and hence the arcs in $\mathcal{H}_s$) are consecutive at $r$. Thus, $\mathcal{H}_s$ is contained in a single strip of $S - \mathcal{H}_q$ and vice versa. Let $\beta, \beta'$ be the arcs in $\mathcal{H}_q$ bounding the unique strip of $S - \mathcal{H}_q$ containing $\mathcal{H}_s$, and let $\gamma, \gamma'$ be the arcs in $\mathcal{H}_s$ bounding the unique strip of $S - \mathcal{H}_s$ containing $\mathcal{H}_q$ (note that we do not exclude the possibility that $\beta = \beta'$ or $\gamma = \gamma'$). Then the complement of $\mathcal{H}$ consists of the remaining strips of $S - \mathcal{H}_q$, $S - \mathcal{H}_s$ and a square bounded by $\beta, \beta', \gamma, \gamma'$. 

\item \paragraph{\textbf{Step 8. $|\mathcal{H}| \leq n-1$.}} Indeed, $|\mathcal{H}|$ is $2$ larger than the number of strips of $S - \mathcal{H}$, so it suffices to show that there are at most $n-3$ of these strips. This is true by Step 7 since $S$ has area $2\pi (n-2)$, and a punctured strip and a square each have area at least $2\pi$. 

\item \paragraph{\textbf{Step 9. $|\mathcal{I}'| \leq n-2$.}} To show this, we intersect $\mathcal{H}$ with a small circle $C$ centered at $r$. Each element of $\mathcal{I}'$ is determined by a pair of points of this intersection, and we connect them by a straight line segment. We also draw a line segment between the outermost points on $C$ corresponding to elements of $\mathcal{H}_q$. By Step 6, these line segments are pairwise intersecting, so by Lemma~\ref{erdos}, $|\mathcal{I}'|+1 \leq |\mathcal{H}| \leq n-1$.
\end{proof}

\section{Proof of Theorem~\ref{maintheorem}}

\begin{figure}
\scalebox{.5}{
\includegraphics[trim=0 520 0 0, clip]{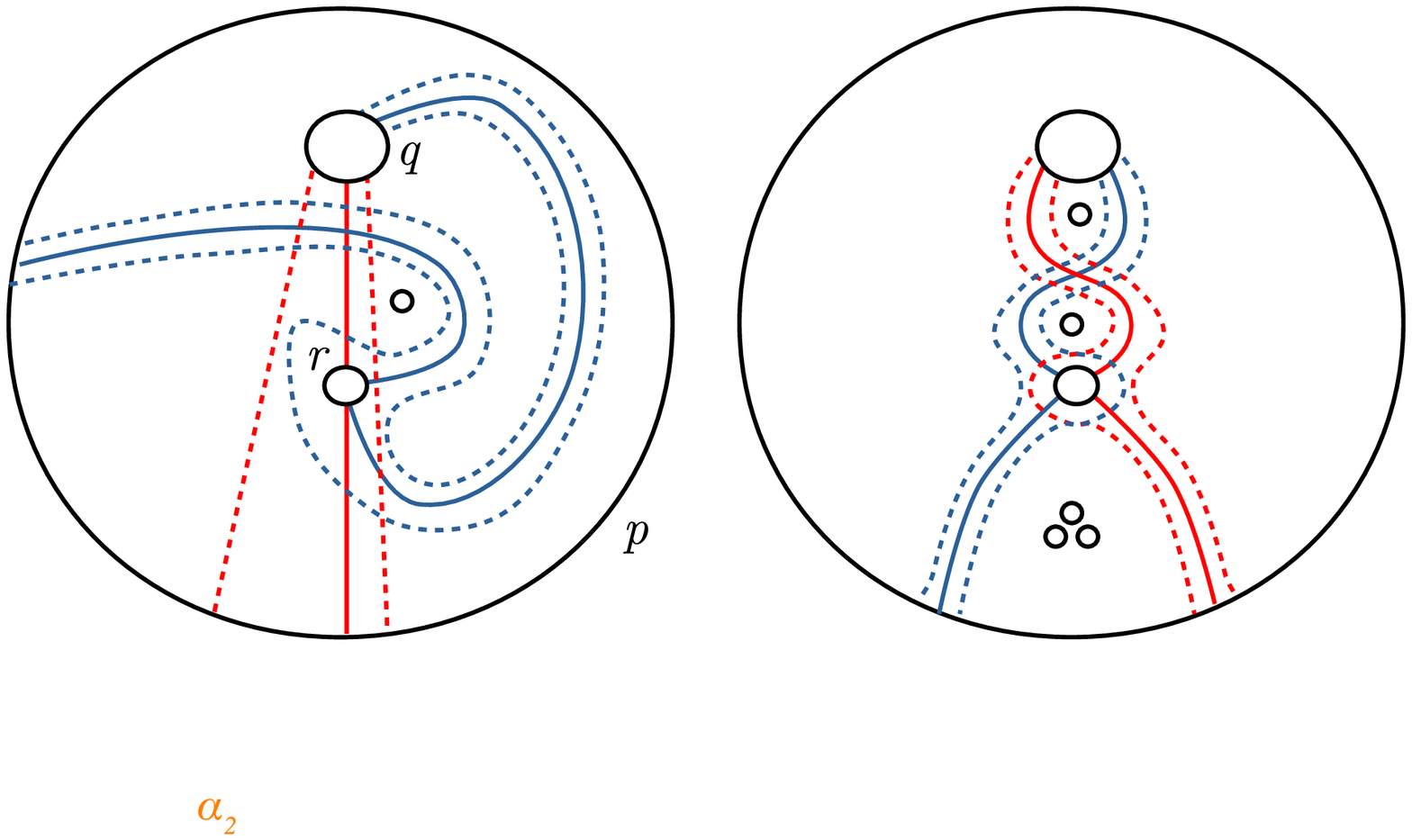}
}
\captionof{figure}{}
\label{fig:alternating}
\end{figure}

\begin{proof}[Proof of Theorem~\ref{maintheorem}]
We proceed by induction on $n$. The case $n=3$ is trivial. Now let $n \geq 4$, and assume the theorem holds if $S$ has fewer punctures. Let $r$ be a puncture of $S$ distinct from $p,q$. Let $\bar{S}$ be the $(n-1)$-punctured sphere obtained from $S$ by forgetting $r$. For each arc $\alpha \in \mathcal{A}$, let $\bar{\alpha}$ be the corresponding arc on $\bar{S}$, and let $\bar{\mathcal{A}} = \{\bar{\alpha} \> : \> \alpha \in \mathcal{A}\}$. We tighten the arcs of $\bar{\mathcal{A}}$ to geodesics. Note that $\bar{\mathcal{A}}$ is a $2$-system on $\bar{S}$, and so $|\bar{\mathcal{A}}| \leq \binom{n-1}{3}$ by the induction hypothesis. Thus, it is enough to show that $|\mathcal{A}| - |\bar{\mathcal{A}}| \leq \binom{n-1}{2}$. To that end, we examine the extent to which the map $\pi: \mathcal{A} \rightarrow \bar{\mathcal{A}}$, $\alpha \mapsto \bar{\alpha}$ is injective. 

By Lemmas~\ref{rhomotopic}, \ref{rhomotopicbigon}, and \ref{extend}, we may add arcs to $\mathcal{A}$ so that for each $\alpha \in \mathcal{A}$,
\[
|\pi^{-1}(\bar{\alpha})| -1 = |\{ \{\alpha_1, \alpha_2\} \in \pi^{-1}(\bar{\alpha}) \> : \> \alpha_1, \alpha_2 \text{ distinct  and disjoint} \}|
\]
Let $\mathcal{P}$ (resp., $\mathcal{Q}$) be the collection of all geodesic arcs $\alpha$ on $S$ starting at $r$ and ending at $p$ (resp., ending at $q$) such that $\alpha$ is contained entirely in a strip bounded by a pair of distinct, disjoint $r$-homotopic arcs in $\mathcal{A}$. Let $\mathcal{R} \subset \mathcal{P} \times \mathcal{Q}$ be the relation consisting of all pairs $(\alpha, \beta)$ such that both $\alpha$ and $\beta$ lie in a single such strip. We claim that $\mathcal{P}, \mathcal{Q}, \mathcal{R}$ satisfy the conditions of Lemma~\ref{relation}, so that $|\mathcal{R}| \leq \binom{n-1}{2}$. Since $|\mathcal{R}| = |\mathcal{A}| - |\bar{\mathcal{A}}|$, this completes the proof.

We first note that for $(\alpha, \beta), (\alpha', \beta') \in \mathcal{R}$ corresponding to pairs of disjoint $r$-homotopic arcs $\gamma_1, \gamma_2 \in \mathcal{A}$ and $\gamma_1', \gamma_2' \in \mathcal{A}$, respectively, we have for some $i, j \in \{1,2\}$ that $r$ produces at least one point of intersection between $\gamma_i, \gamma_j'$. Each point of intersection between the arcs $\alpha, \beta, \alpha', \beta'$ produces an additional point of intersection between $\gamma_i, \gamma_j'$. It follows that there is at most one point of intersection between any two pairs of arcs in $\mathcal{R}$. 

We now show that no arc in $\mathcal{P}$ intersects an arc in $\mathcal{Q}$. Indeed, suppose we have $(\alpha, \beta), (\alpha', \beta') \in \mathcal{R}$ such that $\alpha$ intersects $\beta'$. By the above, $\alpha$ intersects $\beta'$ exactly once and that is the only point of intersection between the pairs of arcs $(\alpha, \beta), (\alpha', \beta')$. But then we can find two arcs in $\mathcal{A}$ that intersect thrice, as shown in Figure~\ref{fig:alternating} (left).

Finally, if there are two intersecting pairs of arcs in $\mathcal{R}$ whose cyclic order around $r$ is not alternating, then we can also find two arcs in $\mathcal{A}$ that intersect thrice, as shown in Figure~\ref{fig:alternating} (right).
\end{proof}

\bibliographystyle{amsalpha}
\bibliography{biblio}

\end{document}